\input amstex
\input xy
\xyoption{all}
\documentstyle{amsppt}
\magnification 1100
\pagewidth{450pt}

\define\Ph{\varphi}
\define\G{\text{\rm Gal}}\define\Dc{\bold D_{\text{\rm cris}}}

\define\res{\text{\rm res}}

\define \boB{\bold B}

\define \A{\bold A}
\define \E{\bold E}
\define \bD{\bold D}

\define \cl{\text{\rm cl}}

\define \Iw{\text{\rm Iw}}

\define\ep{\varepsilon}

\define\Bc{\bold B_{\text{\rm cris}}}

\define\Hom{\text {\rm Hom}}
\define\Ext{\text {\rm Ext}}\define\M {\Cal M}

\define\F{\text{\rm Fil}}

\define\Tam{\text{\rm Tam}}

\define\Ddagrig{\bold D^{\dagger}_{\text{\rm rig}}}
\define\Bdagk{\bold B_{K}^{\dagger}}
\define\Brigdag{\Cal R(K)}


\define\iso{\overset{\sim}\to{\rightarrow}}

\define\Gal{\text{\rm Gal}}

\define \Dst{\bold D_{\text{\rm st}}}
\define \CDcris{\Cal D_{\text{\rm cris}}}
\define  \CDst{\Cal D_{\text{\rm st}}}
\define \CDpst{\Cal D_{\text{\rm pst}}}
\define  \CDdr{\Cal D_{\text{\rm dR}}}
\define \gr{\text{\rm gr}}

\define \CR{\Cal R}
\define \sha{\hbox{\text{\bf III}\hskip -11pt\vrule width10pt
depth0pt height0.4pt\hskip 1pt}}
\topmatter
\title
A generalization of Greenberg's $\Cal L$-invariant
\endtitle
\author
{Denis Benois
\flushpar
$$
$$
}
\endauthor
\abstract Using the theory of $(\Ph,\Gamma)$-modules we generalize Greenberg's construction of the $\Cal L$-invariant 
to  $p$-adic representations which are semistable at $p$.
\endabstract
\leftheadtext{}
\address
Math\'ematiques et Informatique, Universit\'e
Bordeaux I, 351 cours de la Lib\'eration, 33405,
\linebreak
Talence Cedex, France 
\endaddress
\date April 2009
\enddate
\email denis.benois\@math.u-bordeaux1.fr
\endemail
\toc \nofrills {\bf Table of contents}
\widestnumber \head{\S 5.} 
\widestnumber\subhead{5.2}
\specialhead{Introduction} \endspecialhead

\head \S 1. Cohomology of $(\Ph,\Gamma)$-modules
\endhead
\subhead 1.1. Preliminaries
\endsubhead
\subhead 1.2. Semistable and crystalline $(\Ph,\Gamma)$-modules
\endsubhead
\subhead 1.3. Triangulation of $(\Ph,\Gamma)$-modules
\endsubhead
\subhead 1.4. Crystalline and semistable extensions
\endsubhead
\subhead 1.5. Semistable modules of rank $1$
\endsubhead
\head \S 2. The $\Cal L$-invariant
\endhead
\subhead 2.1. Definition of the $\Cal L$-invariant
\endsubhead
\subhead 2.2. $p$-adic $L$-functions  
\endsubhead
\endtoc
\nologo
\endtopmatter
\document

\head
{\bf Introduction}
\endhead

{\bf 0.1.} Let $f=\underset{n=1}\to{\overset \infty \to \sum} a_nq^n$ be a normalized newform of weight $2k$ on $\Gamma_0 (Np)$ where $(N,p)=1$ and let 
$L(f,s)=\underset{n=1}\to{\overset \infty \to \sum}a_nn^{-s}$  denote the complex $L$-function associated to $f$. The Euler factor  of $L(f,s)$ at $p$ can be written in the form $E_p(f,p^{-s})$ where
$E_p(f,X)=\underset i\to {\overset d \to\prod} (1-\alpha_iX)$ with $d\leqslant 2.$ To any $\alpha \in \{\alpha_i \vert 1\leqslant i\leqslant d\}$
such that $v_p(\alpha)<2k-1$ one can associate a $p$-adic $L$-function $L_{p,\alpha}(f,s)$  interpolating $p$-adically the   special values  
$L(f,j)/\Omega_f$ ($1\leqslant j\leqslant 2k-1)$ where $\Omega_f$ denotes the Shimura period of $f$. 

Assume that $U_p(f)=p^{k-1}f$ where $U_p$ is the Atkin-Lehner operator. Then  $E_p(f,X)=1-p^{k-1}X$ and we denote by 
$L_p(f,s)$ the $p$-adic $L$-function associated to
$\alpha=p^{k-1}.$ The interpolation property forces  $L_p(f,s)$ to vanish at $s=k.$  In \cite{MTT} Mazur, Tate and Teitelbaum conjectured that 
there exists an invariant $\Cal L(f)$ which depends only on the restriction of the $p$-adic Galois representation $V_f$ attached to $f$ to a decomposition group at $p$ and such that
$$
L'_{p}(f,k)=\Cal L(f)\,\frac{L(f,k)}{\Omega_f}. \tag{1}
$$
In the weight two case $f$ corresponds to an elliptic curve $E/\Bbb Q$ having split multiplicative reduction at $p$.
The  $p$-adic representation $V_f$ is ordinary  and seats in an exact sequence
of the form
$$
0@>>>\Bbb Q_p(1)@>>>V_f@>>>\Bbb Q_p@>>>0.
$$ 
The class of this extension in $H^1(\Bbb Q_p,\Bbb Q_p(1))$ coincides with the image of the Tate invariant $q_E$ under 
the Kummer homomorphism
 and Mazur, Tate and Teitelbaum conjectured that $\Cal L(f)=\dsize\frac{\log_p q_E}{\text{ord}_p q_E}.$ 
In the higher weight case several  definitions were proposed \cite{Tm}, \cite{Co}, \cite{Mr1}, \cite{O}, \cite{Br}.
It is now known that all these invariants  are equal.  Remark that the Fontaine-Mazur $\Cal L$-invariant  \cite{Mr1} is defined in terms of the filtered $(\Ph,N)$-module $\Dst (V_f)$ and has a natural interpretation in the theory of $(\Ph,\Gamma)$-modules \cite{Cz4}. The  conjecture (1)  was first proved by Greenberg and Stevens \cite{GS} in the  weight two case. In \cite{S}, Stevens generalized this proof
to the higher weights. Other  proofs were found by Kato, Kurihara and Tsuji, Orton, Emerton, ...   
and we refer to \cite{Cz3} and \cite{BDI} for further information and references.
On the other hand,in \cite{G} Greenberg  defined an  
$\Cal L$-invariant  for  pseudo-geometric representations which are { ordinary} at $p$ and suggested a natural generalization of the conjecture(1). Important results in this direction were 
recently proved by Hida \cite{Hi}. 

{\bf 0.2.} In this paper we propose  a definition of the  $\Cal L$-invariant for
representations which are semistable at $p$ generalizing the both Fontaine's and Greenberg's
constructions. In  the subsequent paper \cite{Ben2}   is proved that  this definition
is compatible with Perrin-Riou's theory of $p$-adic $L$-functions \cite{PR}.
The main technical tool is the theory of $(\Ph,\Gamma)$-modules \cite{F1}, \cite{Cz2}. 
We make use of Colmez's  observation that the $(\Ph,\Gamma)$-module 
associated to an irreducible $p$-adic representation may be reducible 
in the category of  $(\Ph,\Gamma)$-modules over the Robba ring $\Cal R$ \cite{Cz4}. 
In particular, semistable  representations are trianguline \cite{BC} and 
we follow Greenberg's approach using the  cohomology of $(\Ph,\Gamma)$-modules
instead Galois cohomology. In \S1,  
for any $(\Ph,\Gamma)$-module $D$ over $\CR$ we define a subgroup $H^1_f(D)$ of $H^1(D)$ which generalizes
$H^1_f(\Bbb Q_p,V)$ of Bloch and Kato \cite{BK} and transpose some classical properties of these groups
to our situation. The proofs are not difficult and follow from fundamental results of Berger \cite{Ber1}, \cite{Ber2},
but do not  seem to be in the literature and we give them in all details. In \S2 the  $\Cal L$-invariant is defined, 
and some related results and conjectures are discussed.

{\bf Acknowledgements.} The author is very grateful to J. Nekov\'a\v r 
for a number of stimulating discussions concerning this project.
 A part of it was done  during the author's visit to the Concordia University at Montreal.
The author would like to thank  A. Iovita for his hospitality and several
interesting discussions.

\head
{\bf \S1. Cohomology  of $(\Ph,\Gamma)$-modules}
\endhead

{\bf 1.1. Preliminaries.}

{\bf 1.1.1. Rings of $p$-adic periods} (see \cite{F2}, \cite{Ber1}, \cite{Cz2}). 
Let $K$ be a finite extension of $\Bbb Q_p .$ We write $K_0$ for the maximal
unramified subextension of $K$, $O_K$ for  the ring of integers of $K$  and  $k_K$ for its residue field.
Let $\sigma $ denote the absolute Frobenius of $K_0/\Bbb Q_p$.
 For any perfect ring $A$ we write  $W(A)$ for the ring of Witt vectors with
 coefficients in $A$.  In particular  $O_{K_0}=W(k_K).$ 
Fix an algebraic closure   $\overline K/K$  and set $G_K=
\text{\rm Gal}(\overline K/K).$  We denote by 
 $C$   the  completion of $\overline K$ and write  
$v_p\,\,:\,\,C@>>>\Bbb R\cup\{\infty\}$ for the $p$-adic valuation normalized so that $v_p(p)=1$. Set
$\vert x\vert_p=\left (\frac{1}{p}\right )^{v_p(x)}.$ 
Let $\mu_{p^n}$  denote the group of $p^n$-th roots of  unity.
Fix a system of primitive roots of  unity  $\ep=(\zeta_{p^n})_{n\geqslant 0}$,
 such that $\zeta_{p^n}^p=\zeta_{p^{n-1}}$ for all $n \geqslant 1$.
Set $K_n= K (\zeta_{p^n})$, $K_{\infty}= \bigcup_{n=0}^{\infty}K_n$, $H_K=\Gal (\bar K/K_\infty)$, $\Gamma =\G(K_{\infty}/K)$
and denote by $\chi \,:\,\Gamma @>>>\Bbb Z_p^*$ the cyclotomic character.

Consider the projective limit 
$$
\widetilde {\bold E}^+=\varprojlim_{x\mapsto x^p} O_C/\,p\,O_C\,
$$
where $O_C$ is the ring of integers of $C.$
Then $\widetilde {\bold E}^+$ is the set of  $x=(x_0,x_1,\ldots ,x_n,\ldots )$ such that  $x_n\in O_C/\,p\,O_C$ and
$x_{n+1}^p=x_n$ for all $n$. For each $n$ choose a lifting $\hat x_n\in O_C$ of $x_n$.
Then for all $m\geqslant 0$ the sequence $\hat x_{m+n}^{p^n}$ converges to 
$x^{(m)}=\lim_{n\to \infty} \hat x_{m+n}^{p^n}\in O_C$
which does not depend on the choice of  $\hat x_n.$
 The ring $\widetilde {\bold E}^+$ equipped with the valuation $v_{\bold E}(x)=v_p(x^{(0)})$ is a
complete local ring of characteristic $p$ with  residue field
 $\overline{\Bbb F}_p$. Moreover it is integrally closed in his field
of fractions $\widetilde {\bold E}=\text {\rm Fr}(\widetilde {\bold E}^+)$. 

Let $\widetilde \A=W(\widetilde \E)$ be the ring of Witt vectors with coefficients
in $\widetilde \E$. Denote by $[\,\,]\,:\,\widetilde \E@>>>W(\widetilde \E)$  the Teichmuller lift.
Any $u=(u_0,u_1,\ldots )\in \widetilde \A$ can be written in the form
$$
u=\underset{n=0}\to{\overset{\infty}\to \sum} [u^{p^{-n}}]p^n.
$$

Set $\pi=[\ep]-1$, $\A_{K_0}^+=O_{K_0} [[\pi]]$ and denote by $\A_{K_0}$  the $p$-adic completion
of $\A_{K_0}^+\left [1/{\pi}\right ]$. This is
a discrete valuation ring with residue field $\E_{K_0}=k_K((\ep-1))$. 
Let $\widetilde\boB= \widetilde \A\left [{1}/{p}\right ]$ and  $\boB_{K_0}=\A_{K_0}\left [{1}/{p}\right ]$ and
let $\boB$ denote  the completion of the maximal unramified extension of $\boB_{K_0}$ in $\widetilde \boB$.
Set $\A=\boB\cap \widetilde \A$, $\widetilde \A^+=W(\E^+)$, $\A^+= \widetilde \A^+\cap \A$ and $\boB^+=\A^+\left [{1}/{p}\right ].$
All these rings are endowed with natural  actions of the Galois group $G_K$ and   Frobenius $\Ph$.

Set $\A_K=\A^{H_K}$ and  $\boB_K=\A_K\left [{1}/{p}\right ].$
Remark that $\Gamma$ and $\Ph$ act on $\boB_{K_0}$ by
$$
\aligned
& \tau (\pi)=(1+\pi)^{\chi (\tau)}-1,\qquad \tau \in \Gamma\\
&\Ph (\pi)=(1+\pi)^p-1.
\endaligned
$$
For any $r>0$ define
$$
\widetilde {\bold B}^{\dagger,r}\,=\,\left \{ x\in \widetilde
{\bold B}\,\,|\,\, \lim_{k\to +\infty} \left (
v_{\E}(x_k)\,+\,\dsize \frac{pr}{p-1}\,k\right )\,=\,+\infty
\right \}.
$$
Set ${\bold B}^{\dagger,r}=\boB \cap \widetilde\boB^{\dagger,r}$,
$\boB_{K}^{\dagger,r}=\boB_{K} \cap \boB^{\dagger,r}$, 
${\bold B}^{\dagger}=\underset{r>0}\to \cup \boB^{\dagger,r}$
and $\bold B^\dag_K=\underset{r>0}\to \cup \boB_{K}^{\dagger,r}$.

It can be shown that
$$
 \boB_{K_0}^{\dagger,r}=\left \{ f(\pi)=\sum_{k\in \Bbb Z}
a_k\pi^k\,\mid \, \text{\rm $a_k\in K_0$ and $f$ is holomorphic and bounded on $p^{-1/r}\leqslant |X|_p<1$} \right \}.
$$
More generally, let $F$ be the maximal unramified subextension of $K_\infty/K_0$ and let
$e= [K_\infty : K_0(\zeta_{p^\infty})].$ Then there exists $r(K)$ and $\pi_K\in \boB^{\dag, r(K)}_K$
such that for any $r\geqslant r(K)$ one has
$$
 \boB_{K}^{\dagger,r}=\left \{ f=\sum_{k\in \Bbb Z}
a_k\pi_K^k\,\mid \, \text{\rm $a_k\in F$ and $f$ is holomorphic and bounded on $p^{-1/er}\leqslant |X|_p<1$} \right \} .
$$ 
Define
$$
\boB^{\dag,r}_{\text{rig},K}\,=\,\left \{ f=\sum_{k\in \Bbb Z}
a_k\pi_K^k\,\mid \, \text{\rm $a_k\in F$ and $f$ is holomorphic  on $p^{-1/er}\leqslant |X|_p<1$} \right \}
\quad \text{\rm if $r\geqslant r(K)$} 
$$ 
and set $\Brigdag =\underset{r\geqslant r(K)}\to \cup \boB^{\dag,r}_{\text{rig},K}.$ 
It is not difficult to
show that $\boB^{\dag}_K$ and $\Brigdag $  are stable under the actions of $\Gamma$ and $\Ph.$
\newline
\,

{\bf 1.1.2. $(\Ph,\Gamma)$-modules} (see \cite{F1}, \cite{Cz2}, \cite{Cz4}). 
 Let $A$ be a commutative ring equipped with a Frobenius  $\Ph .$ A finitely generated free $A$-module
$M$ is said to be a $\Ph$-module if it is equipped with a $\Ph$-semilinear map $\Ph \,:\,M@>>>M$
such that the induced map $A \otimes_{\Ph} M@>>>M$ is an isomorphism. 
If $r\in \Bbb Q$ is a rational written in the form   $r=a/b$ such that $a\in \Bbb Z$, $b\in \Bbb N^*$ and $(a,b)=1$,
we denote by $D^{[r]}$ the $\Ph$-module of rank $b$ defined  by
$$
\Ph (e_1)=e_2,\,\,\Ph (e_2)=e_3, \ldots , \Ph (e_{b-1})=e_b,\,\,\Ph (e_b)=p^ae_1.
$$
A $\Ph$ module is said to be elementary, if it is isomorphic to $D^{[r]}$ for some $r.$

Let $\bold k$ be an algebraically closed field and let
$ \Cal K= W(\bold k)\left [1/p\right ].$ The theory of $\Ph$-modules over $\Cal K$ goes back to Dieudonn\'e and Manin \cite{Mn}. 
Namely, for any $\Ph$-module $D$ over $\Cal K$ there exists a unique decomposition into a direct sum
of elementary modules
$$
D\simeq  \underset{i\in I}\to \oplus D^{[r_i]}.
$$ 
The rationals $r_i$ are called  slopes of $D$. One says that $D$ is pure of slope $r$
if all $r_i=r.$ In particular, if $D$ is a $\Ph$-module over $\Bdagk$, then it can be decomposed
 over $\widetilde {\bold B}$:
$$
D\otimes_{\Bdagk}\widetilde \boB \,\simeq \,\underset{i\in I}\to \oplus D^{[r_i]}.
$$
We say that a $\Ph$-module over $\Bdagk$ is etale if it is pure of slope $0$. 

The analogue of this theory  over the Robba ring $\Brigdag$ is higher non trivial. It was found by Kedlaya \cite{Ke}.
Let $D$ be a $\Ph$-module over $\Brigdag$. Then there exists  a canonical filtration
$$
D_0\subset D_1\subset \cdots \subset D_h=D
$$
and for any $1\leqslant i \leqslant h$, a unique
$\Bdagk$-submodule $\Delta_i (D)$ of $D_i/D_{i-1}$
satisfying the following properties:

i) $\Delta_i(D)$ is a pure  $\Ph$-module of slope $r_i$ such that
$D_i/D_{i-1}=\Delta_i (D)\otimes_{\Bdagk}\Brigdag$.

ii) One has
$$
r_1<r_2<\cdots <r_h.
$$
We say that  $D$  is pure of slope $r$ if $h=1$ and $r_1=r$.  In this case  there exists a unique pure $\Bdagk$-module $\Delta (D)$ of slope $r$ such that $D=\Delta
(D)\otimes_{\Bdagk}\Brigdag .$  

Now assume that $A$ is  a commutative ring  endowed  with  actions of $\Ph$
and $\Gamma$ commuting to each other. A $(\Ph,\Gamma)$-module over $A$ is a
$\Ph$-module equipped with a semilinear action of $\Gamma$ commuting with $\Ph .$

\proclaim{Proposition 1.1.3 } The functors $\Delta \mapsto \Delta
\otimes_{\Bdagk}\Brigdag $ and $D\mapsto \Delta (D)$ are
quasi-inverse equivalences between the category of etale
$(\Ph,\Gamma)$-modules over $\Bdagk$ and the category of
$(\Ph,\Gamma)$-modules over $\Brigdag$ of slope $0$.
\endproclaim
\demo{Proof} see \cite{Cz4}, Proposition 1.4 and Corollary 1.5.
\enddemo

\proclaim{ Proposition 1.1.4} i)  The functor 
$$
\bD^{\dagger}\,\,:\,\,V \mapsto \bD^{\dagger}(V)=(\boB^{\dagger}\otimes_{\Bbb Q_p}V)^{H_K}
$$
establishes an equivalence between the category of $p$-adic representations
of $G_K$ and the category of etale $(\Ph,\Gamma)$-modules over $\Bdagk .$

ii) The functor $\Ddagrig (V)=\Brigdag \otimes
_{\Bdagk} \bD^{\dagger}(V)$ 
establishes an equivalence between the
category of $p$-adic representations of $G_K$ and the category of
$(\Ph,\Gamma)$-modules over $\Brigdag$ of slope $0$.
\endproclaim
\demo{Proof} The first statement  is Fontaine's classification of  $p$-adic representations
\cite{F1} together with the main theorem of  \cite{CC}. The second statement follows
from i) together with proposition 1.1.3. See \cite{Cz4}, Proposition 1.7 for details.
\enddemo

{\bf 1.1.5. Cohomology of $(\Ph,\Gamma)$-modules} (see \cite{H1}, \cite{H2}, \cite{Li}).  Fix a generator $\gamma$ of $\Gamma$. If $D$ is a $(\Ph,\Gamma)$-module over a ring $A$, we denote by
$ C_{\Ph,\gamma}(D)$ the complex
$$
C_{\Ph,\gamma} (D)\,\,:\,\,0@>>>D @>f>> D\oplus D@>g>> D@>>>0
$$
where $f(x)=((\Ph-1)\,x,(\gamma -1)\,x)$ and
$g(y,z)=(\gamma-1)\,y-(\Ph-1)\,z.$ We shall write  $H^i(D)$  for the cohomology
of $C_{\Ph,\gamma}(D).$ A short
exact sequence of $(\Ph,\Gamma)$-modules
$$
0@>>>D'@>>>D@>>>D''@>>>0
$$
gives rise to a long exact sequence of cohomology groups:
$$
0@>>>H^0(D')@>>>H^0(D)@>>>H^0(D'')@>>>H^1(D')@>>>\cdots @>>>
H^2(D'')@>>>0.
$$
If $D_1$ and $D_2$ are  $(\Ph,\Gamma)$-modules define a bilinear map
$$
\cup \,\,:\,\,H^i(D_1)\times H^j(D_2) @>>>H^{i+j}(D_1\otimes D_2)
$$
by
$$
\align
&\cl (x_1) \cup \cl (x_2) =\cl (x_1\otimes x_2)\qquad \text{if $i=j=0$},\\
&\cl (x_1) \cup \cl (x_2,y_2)=\cl (x_1\otimes x_2,x_1\otimes y_2) \qquad \text{if $i=0$, $j=1$},\\
&\cl (x_1,y_1)\cup \cl (x_2,y_2)=\cl (y_1\otimes \gamma (x_2)-x_1\otimes \Ph (y_2)) \qquad \text{if $i=1$, $j=1$}.
\endalign
$$
For any $(\Ph,\Gamma)$-module $D$ let  $D(\chi)$ denote the $\Ph$-module $D$ endowed with the 
action of $\Gamma$  on $D$ twisted by the cyclotomic character $\chi .$  
Set $D^*=\Hom_{\Brigdag} (D,\Brigdag).$

The following theorem extends  the main results of  \cite{H1}, \cite{H2} 
to  $(\Ph,\Gamma)$-modules over $\Brigdag .$ 
 
\proclaim{Theorem 1.1.6} If $D$ be a $(\Ph,\Gamma)$-module over $\Brigdag$ 
then

i)  $H^i(D)$ are finite dimensional $\Bbb Q_p$-vector spaces and 
$$
\chi (D)=\underset{i=0}\to{\overset 2\to \sum} (-1)^i \dim_{\Bbb Q_p} H^i(D)\,=\,-[K:\Bbb Q_p]\,\text{\rm rg} (D)
$$

ii) For $i=0,1,2$ the cup product
$$
H^i(D)\times H^{2-i}(D^*(D(\chi))) @>\cup >> H^2(\Brigdag (\chi))
$$
is a perfect pairing into $H^2(\Brigdag (\chi))\simeq \Bbb Q_p.$

iii) If  $D$ is an etale $(\Ph,\Gamma)$-module over
$\Bdagk $ then the natural map
$$
C_{\Ph,\gamma}(D) @>>> C_{\Ph,\gamma}(D\otimes_{\Bdagk} \Brigdag)
$$
is a quasi isomorphism. In particular, if $V$ is p-adic representation of $G_{K}$,  the Galois cohomology $H^*(K,V)$ is canonically and
functorially isomorphic to $H^*(\Ddagrig (V)).$
\endproclaim   
\demo{Proof}  See  \cite{Li}, Theorems 1.1 and 1.2. 
\enddemo

{\bf Remark 1.1.7.} The isomorphism $H^2(\CR (K)(\chi)) \simeq \Bbb Q_p$ is  constructed  in \cite{H2}. 
If $K$ is  unramified over $\Bbb Q_p$ it is given by the formula
$$
\cl (\alpha) \mapsto -\left(1-\frac{1}{p} \right)^{-1}(\log \chi (\gamma))^{-1} \,\res \frac{\alpha d\pi}{1+\pi}.
$$
It can be shown that this definition is compatible with the canonical isomorphism
$H^2(\Bbb Q_p,\Bbb Q_p(1))\simeq \Bbb Q_p$  of the local class field theory (see \cite{Ben1}, section 2.2). 
\newline
\,

{\bf 1.2. Semistable and crystalline $(\Ph,\Gamma)$-modules} (see \cite{Ber1}, \cite{Ber2}, \cite{BC}).

{\bf 1.2.1.} Let $\log \pi$ be a transcendental element over the field of fractions of $\Brigdag$
equipped with the following actions of $\Ph$ and $\Gamma $ 
$$
\aligned 
&\Ph (\log \pi)=p\log \pi +\log \left (\frac{\Ph (\pi)}{\pi^p}\right )\,,\\
& \gamma (\log \pi)=\log \pi + \log  \left (\frac{\gamma (\pi)}{\pi}\right ).
\endaligned
$$ 
Remark that the series $\log \left ( {\Ph (\pi)}/{\pi^p}\right )$ and $\log  \left ({\gamma (\pi)}/{\pi}\right )$
converge in $\Brigdag .$ Set $\CR_{\log}(K) =\Brigdag [\log \pi]$ and define a monodromy operator 
$N\,:\,\CR_{\log}(K) @>>>\CR_{\log}(K)$ by $N=-\left (1-\dsize\frac{1}{p} \right )^{-1}\dsize\frac{d}{d \log \pi}.$
Write  $K_\infty ((t))$ for  the ring of Laurent power series 
equipped with the filtration $\F^i K_\infty ((t)) =t^iK_\infty [[t]]$ and the
natural  action of $\Gamma$ given by $\gamma \left (\sum a_it^i \right )\,=\,
\sum \gamma (a_i) \chi (\gamma)^i t^i .$  
Recall that for any $n\geqslant 0$ and $r_n=p^{n-1}(p-1)$ there exists a well defined injective 
homomorphism
$$
\iota_{n}=\Ph^{-n} \,\,:\,\, \bold B^{\dag, r_n}_{\text{\rm rig},K}  @>>> K_\infty [[t]]
$$
which is characterized by the fact that  $\iota_{n} (\pi)\,=\,\zeta_{p^n}e^{t/p^n}-1 $ (see for example \cite{Ber1}, \S 2.4).
For any $r>0$ let $n(r)$ denote the smallest integer $n$  such that $r_n\geqslant r.$ 

\proclaim{Lemma 1.2.2} Let $D$ be a $(\Ph,\Gamma)$-module over $\Brigdag$.  There exists $r(D)>0$
such that for any $r \geqslant r(D)$ there exists a unique free $\bold B^{\dag, r}_{\text{\rm rig},K} $-submodule  $D^{(r)}$ of $D$ stable under
$\Gamma$ and  having the following properties:

i) If  $s>r$ then $D^{(s)}\simeq D^{(r)}\otimes_{\bold B^{\dag, r}_{\text{\rm rig},K}}\bold B^{\dag, s}_{\text{\rm rig},K} $.

ii) $D^{(r)}\otimes_{\bold B^{\dag, r}_{\text{\rm rig},K}}\Brigdag \simeq D$;

iii) $D^{(rp)}\simeq \bold B^{\dag, rp}_{\text{\rm rig},K}\otimes_{\bold B^{\dag, r}_{\text{\rm rig},K},\Ph} D^{(r)}.$

\endproclaim
\demo{Proof} This is Theorem 1.3.3 of \cite{Ber2}.
\enddemo

{\bf 1.2.3.} Let $D$ be a $(\Ph,\Gamma)$-module over $\Brigdag$. Following \cite{BC},  define
$$
\CDdr (D)\,=\,(K_\infty ((t))\otimes_{\bold B^{\dag, r}_{\text{\rm rig},K},\,\iota_{n}} D^{(r)} )^{\Gamma},\qquad   r\geqslant r(D),\,\, n\geqslant n(r).
$$
From Lemma 1.2.2 it follows that this definition does not depend on the choice of $r$ and $n .$ 
Since $K_\infty ((t))^{\Gamma}= K$, a  standard argument shows that $\CDdr (D)$ is a
$K$-vector space such that
$$   
\dim_{K} \CDdr (D) \leqslant \text{\rm rg} (D).
$$
Moreover, $\CDdr (D)$ is equipped with the induced filtration
$$
\F^i \CDdr (D)\,=\,(t^i K_\infty [[t]]\otimes_{\bold B^{\dag, r}_{\text{\rm rig},K},\,\iota_{n}} D^{(r)} )^{\Gamma},
\qquad n\geqslant n(r).
$$
The jumps of this filtration (with multiplicities) will be called Hodge-Tate weights of $D$ .
Next, we define
$$
\aligned
&\CDcris (D)=\left (D\otimes_{\Brigdag} \Brigdag [1/t] \right )^{\Gamma}\,=\,(D[1/t])^{\Gamma},\\ 
&\CDst (D)= \left (D\otimes_{\Brigdag}\CR_{\log}(K)[1/t] \right )^{\Gamma}.
\endaligned
$$
Then $\CDcris (D)$ is a finite dimensional $K_0$-vector space equipped with a linear action of $\Ph$.
For any $r\geqslant r(D)$ and $n\geqslant n(r)$ the map $\iota_n$ induces a natural inclusion
$
\CDcris (D) \hookrightarrow \CDdr (D)$ and we set
$$
\F^i\CDcris (D)= \Ph^{n} (\iota_n(\CDcris (D)) \cap \F^i \CDdr (D)).
$$
It is easy to see that $\F^i \CDcris (D)$ is  a filtration on $\CDcris (D)$ which does not depend
on the choice of $n$ and $r$. The same arguments show that $\CDst (D)$ is a finite dimensional $K_0$-vector 
space equipped with natural actions of $\Ph$ and $N$ and a filtration $\F^i \CDst (D)$ induced
from $\CDdr (D).$ From definitions it follows immediately that $\CDcris (D)=\CDst (D)^{N=0}.$ 
Moreover 
$$
\dim_{K_0} \CDcris (D) \leqslant \dim_{K_0} \CDst (D)\leqslant  \dim_{K} \CDdr (D) \leqslant 
\text{\rm rg} (D).
$$

\proclaim{Definition} We say that $D$ is crystalline (resp. semistable, resp. de Rham)
if $\dim_{K_0} \CDcris (D)=\text{\rm rg}_{\Brigdag} (D)$ (resp. $\dim_{K_0} \CDst (D)=\text{\rm rg}_{\Brigdag} (D)$,
resp. $\dim_{K} \CDdr (D)=\text{\rm rg}_{\Brigdag} (D)$).
\endproclaim

This definition is motivated by the following proposition which summarizes the results of Berger 
about the classification of $p$-adic representations in terms of $(\Ph,\Gamma)$-modules. 

\proclaim{Proposition 1.2.4} Let $V$ be a $p$-adic representation of $G_{K}.$ Then
$$
\bold D_* (V)\simeq {\Cal D}_* (\Ddagrig (V)), \qquad \text{\rm where}\,\, *\in \{\text{\rm dR},\text{\rm st},\text{\rm cris}\}.
$$
In particular, $V$ is a de Rham (resp. crystalline, resp. semistable) if and only if 
$\Ddagrig (V)$ is.
\endproclaim
\demo{Proof} see \cite{Ber1}, Theorem 0.2 and Proposition 5.9.
\enddemo

{\bf 1.2.5.} Let $L/K$ be a finite extension and $\Gamma_L=\text{\rm Gal}(L_\infty/L).$
Remark that $\Brigdag \subset \CR (L).$ If
$D$ is a $(\Ph,\Gamma)$-module over $\Brigdag$ 
 we set $D_L= \CR (L)\otimes_{\Brigdag}  D $. It is easy to see
 that $D_L$ has a natural structure of $(\Ph ,\Gamma_L)$-module. We say that $D$ is potentially semistable
if there exists $L/K$ such that $D_L$ is semistable. For any $D$ define
$$
\Cal D_{\text{pst}}(D)=\underset{L/K} \to \varinjlim \Cal D_{\text{st}/L} (D_L).
$$
Then $\Cal D_{\text{pst}}(D)$ is a finite dimensional $K_0^{nr}$-vector space equipped
with natural actions of  $\Ph$, $N$ and a discrete action of $G_K$. It is easy to see that $D$ is potentially semistable
if and only if $\dim_{K_0^{nr}} \Cal D_{\text{pst}}(D) = \text{\rm rg} (D).$

\proclaim{Proposition 1.2.6} Let $D$ be a $(\Ph,\Gamma)$-module. The following
statements are equivalent: 

1) $D$ is potentially semistable.

2) $D$ is a de Rham.
\endproclaim
\demo{Proof} If $D$ is etale, this is the $p$-adic monodromy conjecture stated
by Fontaine. In \cite{Ber1}, Corollary 5.22 Berger deduced it from Crew's conjecture
which was proved independently by Andr\'e, Mebkhout et Kedlaya. Remark that Berger's  arguments
work  if $D$ is not more  etale and give therefore a proof of the Proposition. See especially Lemma 5.13
and Propositions 5.14 and 5.15 of \cite{Ber1}.
\enddemo

{\bf 1.2.7.}  A filtered $(\Ph, N)$-module over $K$ is a finite dimensional $K_0$-vector space $M$
equipped with the following structures:

$\bullet$ an exhaustive decreasing filtration $(\F^i M_K)$ on $M_K=K\otimes_{K_0}M$;

$\bullet$ a $\sigma$-semilinear bijective map $\Ph \,:\,M@>>>M$;

$\bullet$ a $K$-linear nilpotent operator $N\,:\,M@>>>M$  such that
$N \,\Ph =p\,\Ph \,N.$ 

A filtered $(\Ph,N, G_K)$-module over $K$ is a $(\Ph,N)$-module $M$ over $K^{nr}$ equipped with a
 semilinear action of $G_K $ which is discrete on the inertia subgroup $I_K=\text{Gal}(\overline K/K^{nr})$ 
 and an  exhaustive decreasing filtration $(\F^i M_{K^{nr}})$ on $M_{K^{nr}}=K^{nr}\otimes_{K_0^{nr}}M$
stable under the action of $G_K.$ In particular the filtration on $M_{K^{nr}}$ is completely defined
by the filtration $\F^i M_K=M_K\cap \F^i M_{K^{nr}}$ on the $K$-vector space $M_K=(M_{K^{nr}})^{G_K}.$
 
A $K$-linear map $f\,:\,M@>>>M'$ is said to be a morphism of filtered $(\Ph,N, G_K)$-modules
if $f$ commutes with $\Ph$, $N$ and $G_K$ and $f(\F^i M_K)\subset \F^i M_K'$ for all $i\in \Bbb Z.$

The category $\bold M\bold F_K^{\Ph,N, G_K}$  of filtered $(\Ph,N)$-modules is additive, has kernels and 
cokernels but is not abelian. Moreover it is equipped with the tensor product defined by
$$
\aligned
&\F^i (M\otimes M')_K =\underset{j+k=i}\to \sum \F^j M_K\otimes_{K^{nr}} \F^k M_K',\\
&\Ph (m\otimes m')=\Ph (m)\otimes \Ph (m'),\\
&N(m\otimes m')= N(m)\otimes m' + m\otimes N(m').
\endaligned
$$
Denote by $\bold 1$ the vector space $K_0^{nr}$ equipped with the natural actions of $\Ph$ and $G_K$ ,
the trivial action of $N$ and the filtration given by
$$
\F^i \bold 1_K =\cases K &\text{ if $i\leqslant 0$},\\
0 &\text{ if $i>0$}.
\endcases
$$
Then $\bold 1$ is a unit object of $\bold M\bold F_K^{\Ph,N,G_K},$ i.e.
$$
\bold 1\otimes M\simeq M\otimes \bold 1\simeq M
$$
for any $M$.
The functor $M \mapsto K_0^{nr}\otimes_{K_0} M$ is a full embedding of $\bold M\bold F_K^{\Ph,N}$
into $\bold M\bold F_K^{\Ph,N,G_K}$.

A filtered Dieudonn\'e module is an object  $M \in \bold M\bold F_K^{\Ph,N}$ such that
the operator $N$ is trivial on $M.$  Filtered Dieudonn\'e modules
form a full subcategory $\bold M\bold F_K^{\Ph}$ of $\bold M\bold F_K^{\Ph,N}.$  
Denote by $\bold M^{\,\Ph,\Gamma}_{\text{cris},K},$ $\bold M^{\,\Ph,\Gamma}_{\text{st},K}$
and $\bold M^{\,\Ph,\Gamma}_{\text{pst},K}$
the categories of crystalline, semistable and potentially semistable  $(\Ph,\Gamma)$-modules respectively.
For all these categories   we will use the following convention. A sequence
of objects
$$
0@>>> A_1@>>>A@>>>A_2@>>>0
$$
is said to be exact if $A/A_1$ is isomorphic to $A_2$.

\proclaim{Lemma 1.2.8} Let
$$
0@>>>D_1 @>>>D@>>>D_2 @>>>0
\tag{2}
$$
be an exact sequence of $(\Ph,\Gamma)$-modules. If $D$ is a de Rham (resp. crystalline,
resp. semistable) , then $D_1$ and $D_2$ are de Rham (resp. crystalline, resp. semistable)
and  the sequence 
$$
0@>>>\Cal D_* (D_1)@>>>\Cal D_*(D)@>>>\Cal D_*(D_2)@>>>0 , \qquad *\in\{\text{\rm dR},\text{\rm cris},\text{\rm st}\}
$$
is exact in the category of filtered $K$-vector spaces (resp. in $\bold{MF}_K^{\Ph}$, resp. in $\bold{MF}^{\Ph,N}_K$).
\endproclaim
\demo{Proof} The proof is standard. Assume that $D$ is a de Rham module. 
By Lemma 2.2.11 of \cite{BC}, for  $r$ big enough
the sequence
$$
0@>>>D_1^{(r)}@>\alpha>>D^{(r)}@>\beta>>D_2^{(r)}@>>>0
$$
is exact. Tensoring this sequence with $\bold B^{\dag,r}_{\text{rig},K}$ and taking invariants,
one obtains an exact sequence of $K$-vector spaces
$$
0@>>>\CDdr (D_1)@>\alpha>>\CDdr (D)@>\beta>>\CDdr (D_2).
$$
 As $\dim_K \CDdr (D)=\text{\rm rg} (D)$ and 
$\dim_K \CDdr (D_i)\leqslant \text{\rm rg} (D_i)$ ($i=1,2$), this implies that
\linebreak
$\dim_K \CDdr (D_i) = \text{\rm rg} (D_i).$
Thus $D_i$ are de Rham modules. Next, the standard argument involving Artin's trick
shows that for any de Rham $(\varphi,\Gamma)$-module $D$ one has an isomorphism  
$$
D^{(r)}\otimes_{\Brigdag ,\,\iota_n} K_\infty ((t))\simeq \,\CDdr (D)\otimes_{K} K_\infty ((t)).
$$
Thus, for any $k\in \Bbb Z$
$$
D^{(r)}\otimes_{\bold B^{\dag,r}_{\text{rig},K} ,\,\iota_n} t^k K_\infty [[t]]\,=\,\underset{i\in \Bbb Z}\to
\sum \F^{-i}\CDdr (D)\otimes_{K}  t^{k+i} K_\infty [[t]].
$$
Since
$$
H^1(\Gamma, K_\infty t^m)\,=\,\cases 0 &\text{if $m\ne 0$}\\
K &\text{if $m=0$}
\endcases
$$
(see \cite{T}, Proposition 8), we obtain that
$$
H^1\left (\Gamma , D^{(r)}\otimes_{\bold B^{\dag,r}_{\text{rig},K},\,\iota_n} t^k K_\infty [[t]]\right )\,=\,
\F^k \CDdr (D).
$$
The short exact sequence (2) induces a long exact sequence
$$
\multline
0@>>> \F^k \CDdr (D_1)@>\alpha>>\F^k \CDdr (D) @>\beta>>\F^k \CDdr (D_2)@>\delta>>\\
H^1 \left (\Gamma , D_1^{(r)}\otimes_{\bold B^{\dag,r}_{\text{rig},K},\,\iota_n} t^k K_\infty [[t]] \right )@>\alpha_*>>
H^1 \left (\Gamma , D^{(r)}\otimes_{\bold B^{\dag,r}_{\text{rig},K},\,\iota_n} t^k K_\infty [[t]]\right )@>>> \cdots .
\endmultline
$$
Since the diagram
$$
\CD 
H^1\left (\Gamma , D_1^{(r)}\otimes_{\bold B^{\dag,r}_{\text{rig},K},\,\iota_n} t^k K_\infty [[t]]\right ) @>\alpha_* >>
H^1\left (\Gamma , D^{(r)}\otimes_{\bold B^{\dag,r}_{\text{rig},K},\,\iota_n} t^k K_\infty [[t]]\right )\\
@V^{\simeq}VV @V^{\simeq}VV  \\
\F^k \CDdr (D_1) @>\text{\rm id}>> \F^k \CDdr (D)
\endCD
$$
is commutative, the map $\alpha_*$ is injective and we obtain that for all $k$ the sequence 
$$
0@>>> \F^k \CDdr (D_1)@>\alpha>>\F^k \CDdr (D) @>\beta>>\F^k \CDdr (D_2)@>>>0
$$
is exact. Thus, the sequence $0@>>>\CDdr (D_1)@>>>\CDdr (D)@>>>\CDdr (D_2)@>>>0$ is exact in the category
of filtered $K$-vector spaces.
 The case of crystalline (resp. semistable) modules is analogous 
and is omitted here.

\enddemo

\proclaim{Proposition 1.2.9} The functors

$$
\left \{
\aligned
&\Cal D_{\text{pst}} \,\,:\,\,\bold M^{\,\Ph,\Gamma}_{\text{pst},K}  @>>> \bold M\bold F_K^{\Ph,N,G_K},\\
&D\mapsto \Cal D_{\text{pst}} (D),
\endaligned
\right.
$$
$$
\left \{
\aligned
&\CDst \,\,:\,\,\bold M^{\,\Ph,\Gamma}_{\text{st},K}  @>>> \bold M\bold F_K^{\Ph,N},\\
&D\mapsto \CDst (D)
\endaligned
\right.
$$
and
$$
\left \{
\aligned
&\CDcris \,\,:\,\,\bold M^{\,\Ph,\Gamma}_{\text{cris},K}  @>>> \bold M\bold F_K^{\Ph},\\ 
& D\mapsto \CDcris (D)
\endaligned
\right.
$$

are equivalences of categories.

\endproclaim
\demo{Proof} This proposition is  a reformulation of the main result of \cite{Ber2}. 
Let $\nabla$ denote the operator $ \dsize \frac{\log (\gamma)}{\log \chi (\gamma)}$ 
and let $\partial = t^{-1}\nabla $.
Remark that $\partial $ acts as $(1+\pi)\,\dsize \frac{d}{d\pi}$ on $\Brigdag .$
In op. cit. Berger establishes an equivalence $\Cal S$ between the category of $(\Ph,\Gamma)$-modules over $\Brigdag$
such that $\text{Lie}(\Gamma)$ acts locally trivially and the category of filtered $(\Ph,N, G_K)$-modules. 
Namely, if $D$ is such a $(\Ph,\Gamma)$-module, then there exists a unique submodule $A\subset D$ such that
$\partial (A)\subset A$ and $A[1/t]=D[1/t]$ and $\Cal S (D)$ is defined as the space of $L$-solutions
of the differential equation $(A,\partial)$ for $L/K$ big enough.  

Consider the natural  injection
$$
\Cal D_{\text{st}/L} (D)\otimes_{L_0} \CR_{\log}(L) [1/t] \hookrightarrow 
D_L\otimes_{{\CR}(L)} \CR_{\log}(L)  [1/t]
$$
If $D$ is semistable over $L$ set $B=(\Cal D_{\text{pst}/L} (D)\otimes_{L_0}  \CR_{\log}(L) )^{N=0}$. Then $B$
is a $(\Ph,\Gamma_L)$-module over $ \CR (L)$ of rank $\text{rg}_{\Brigdag} (D)$
 which is 
stable under $\partial$. Let $A=B^{H_K}.$ By Hilbert theorem 90 we have $B= \CR (L) 
\otimes_{\Brigdag} A$. Then  $\partial (A)\subset A$, $A[1/t]=D[1/t]$ and 
the space of $L$-solutions of the differential equation $(A,\partial)$ is
$$
(A\otimes_{\Brigdag}\CR_{\log}(L) )^{\partial =0} =\Cal D_{\text{st}/L} (D).
$$
Thus, $\Cal S (D)=\Cal D_{\text{st}/L} (D).$ 

Conversely, let $M$ be a filtered $(\Ph,N, G_K)$-module and  let $L/K$ be a finite extension
such that $I_L$ acts trivially on $M.$ Then $M'= M^{G_L}$ is an $L_0$-lattice of $M$ and we 
set
$$
\bold D = (M'\otimes_{L_0} \CR_{\log}(L) )^{N=0}.
$$
It is easy to see that $\bold D$ is a $(\Ph,\Gamma_L)$-module of rank $\dim_{K_0^{nr}} (M).$
By Theorem II.1.2 of \cite{Ber 2} there exists a unique $(\Ph,\Gamma_L)$-module $\Cal M_L(M)\subset \bold D$ such that
$\Cal M_L(M)[1/t]=\bold D[1/t]$ and
$$
L_n[[t]]\otimes_{\bold B^{\dag,r}_{\text{rig},L},\iota_n}
 \Cal M_L(M)^{(r)} =\F^0 (L_n((t))\otimes_{L_0,\Ph^{-n}}  M') 
\qquad \text{for all $n \geqslant n(r)$.}
$$
Set $\Cal M(M)=\Cal M_L(M)^{H_K}.$
Berger proves that the functor  $M \mapsto \Cal M(M)$ is a quasi-inverse to $\Cal S.$ From the isomorphism
$
\bold D\otimes_{\CR{(L)} } \CR_{\text{log}}(L)\simeq  M\otimes_{L_0}\CR_{\text{log}}(L) 
$
it follows easily that
$$
\Cal D_{\text{st}/L} (\Cal M(M)_L)= (\bold D\otimes_{\CR_L}\CR_{\text{log}}(L))^{\Gamma_L} =M'.
$$
Thus $D$ is a semistable $(\Ph,\Gamma)$-module such that  $\Cal D_{\text{pst}} (D)=M$ and we proved that $\Cal D_{\text{pst}}$
is an equivalence between the category of potentially semistable $(\Ph,\Gamma)$-modules and $\bold M\bold F_K^{\Ph,N,G_K}.$ 
Passing to the subcategory $\bold M^{\,\Ph,\Gamma}_{\text{st},K}$  (resp. $\bold M^{\,\Ph,\Gamma}_{\text{cris},K}$ )
we obtain immediately  that $\CDst$ (resp. $\CDcris$) is an equivalence of the category of semistable (resp. crystalline) modules onto 
$\bold M\bold F_K^{\Ph,N}$  (resp. $\bold M\bold F_K^{\Ph}$). 
The proposition is proved. 
\enddemo
\newpage

{\bf 1.3. Triangulation of $(\Ph,\Gamma)$-modules} (see \cite{Ber2}, \cite{Cz4}, \cite{BC}). 

{\bf 1.3.1.} The results of this section will not be used in the remainder of this paper.
Nevertheless the notion of a trianguline representation is closely related to our definition of the
$\Cal L$-invariant and we review it here. For simplicity we assume that $K=\Bbb Q_p$
and write $\CR$ for $\CR (\Bbb Q_p)$
 but
fix a finite extension $L$ of $\Bbb Q_p$ as the coefficient field.  Equip $L$  with trivial 
actions of $\Ph$ and $\Gamma$ and set
$\CR_{L}=\CR_{\Bbb Q_p} \otimes L.$ Remark that the theory of sections 1.1 and 1.2 
extends without difficulty to $(\Ph,\Gamma)$-modules over $\CR_{L}.$   
Let $D$ be  such a  module.  A triangulation of $D$ is a strictly 
increasing filtration
$$
\{0\}=F_0D\subset \cdots \subset  F_i D\subset F_{i+1}D\subset \cdots \subset F_d D=D
$$
by $(\Ph,\Gamma)$-submodules over $\CR_{L}$ such that

$\bullet $ every $F_iD$ is a saturated submodule of $D$;

$\bullet $ the factor modules $\gr_i (D)= F_iD/F_{i-1}D$ are free of rank $1$.
\flushpar
Triangular modules were first studied in \cite{Cz4}.

Now  assume that $D$ is semistable and that all the eigenvalues of $\Ph \,:\,\CDst (D)@>>>\CDst (D)$
are in $L$. Following Mazur \cite{Mr2}, a refinement of $D$ is a filtration  on $\CDst (D)$
$$
\{0\}= \Cal F_0\CDst (D)\subset \Cal F_1\CDst (D)\subset \cdots \subset \Cal F_d\CDst (D)=\CDst (D)
$$
by $L$-subspaces stable under $\Ph$ and $N$ and such that each factor  $\gr_i \CDst (D)$ is of dimension $1$.
Any refinement fixes an ordering $\alpha_1,\ldots ,\alpha_d$
of eigenvalues of $\Ph$ and an ordering $k_1,\ldots ,k_d$ of Hodge weights of $D$ taken with multiplicities.

\proclaim{Proposition 1.3.2} Let $D$ be a semistable $(\Ph,\Gamma)$-module over $\CR_{L}.$ 

i) The equivalence between the category of semistable modules and the category of filtered
$(\Ph,N)$-modules induces a bijection between the set of triangulations of $D$ and the set of refinements of $D$.

ii) If  $(F_i D)_{i\in \Bbb Z}$ is the triangulation associated to a refinement $(\Cal F_i\CDst (D))_{i\in \Bbb Z}$  then 
for each $i$ the factor $\gr_iD$ is isomorphic to $\CR_{L}(\delta_i)$ where $\delta_i$ 
is defined by $\delta_i(p)=\alpha_i p^{-k_i}$ and $\delta_i(u)=u^{-k_i}$ ($u\in \Bbb Z_p^*$).
\endproclaim
\demo{Proof} For crystalline representations this Proposition was proved  in \cite{BC}, Proposition 2.4.1 and the same
proof works in the general case. On the other hand, it can be deduced easily from 
Proposition 1.2.9.  Indeed, the first statement is obvious. Next, let $e_{i}$ denote the canonical generator
of $\CR_{L}(\delta_i)$. Then  $\CDst (\CR_{L}(\delta_i))$
is the one-dimensional $L$ vector space generated by $m_i=t^{k_i}\otimes e_i.$ 
One has $\Ph  (m_i)=\alpha_i m_i$, $N m_i=0$ and $k_i$ is the unique Hodge number of $\CDst (\CR_{L}(\delta_i)).$
Thus $\CDst (\CR_{L}(\delta_i))\simeq \CDst (\gr_i (D))$ and  
$\gr_iD$ is isomorphic to $\CR_{L}(\delta_i)$ by Proposition 1.2.9.
\enddemo

{\bf 1.4. Crystalline and semistable extensions.}

{\bf 1.4.1.} Let $D$ be a $(\Ph,\Gamma)$-module over $\Brigdag .$  As usually, $H^1(D)$ can be described 
in terms of extensions.  Namely, to any cocycle $\alpha =(a,b)\in Z^1(C_{\Ph,\gamma}(D))$ 
we can associate the extension
$$
0@>>>D@>>>D_\alpha @>>>\Brigdag @>>>0
$$
defined by
$$
D_{\alpha}= D\oplus \Brigdag \,e, \qquad (\Ph-1)\,e=a, \quad (\gamma-1)\, e=b\,.
$$ 
This construction gives rise to an isomorphism
$$
H^1 (D) \simeq \text{\rm Ext}^1_{\CR} (\Brigdag ,D).
$$

\proclaim{Definition} We say that $\cl (\alpha)\in H^1(D)$ is crystalline (resp. semistable) if 
$$
{\text{\rm dim}}_{K_0} \CDcris (D_\alpha)={\text{\rm dim}}_{K_0} \CDcris (D)+1
$$
(resp. if ${\text{\rm dim}}_{K_0} \CDst (D_\alpha)={\text{\rm dim}}_{K_0} \CDst (D)+1$) 
and define
$$
\align
&H^1_{f} (D)\,=\,\{ \cl (\alpha) \in H^1(D) \mid \text{\rm $\cl(\alpha)$ is crystalline}\},\\
&H^1_{st} (D)\,=\,\{ \cl (\alpha) \in H^1(D) \mid \text{\rm $\cl(\alpha)$ is semistable}\}.
\endalign
$$

\endproclaim
It is easy to see that $H^1_f(D)\subset H^1_{st}(D)$ are $\Bbb Q_p$-subspaces of $H^1(D).$

\proclaim{Proposition 1.4.2} Let $V$ be a $p$-adic representation of $G_K$. Following
Bloch and Kato \cite{BK} define
$$
\align
&H^1_f(K,V)=\ker (H^1(K,V) @>>>H^1(K,V\otimes_{\Bbb Q_p} \Bc)),\\
&H^1_{st}(K,V)=\ker (H^1(K,V) @>>>H^1(K,V\otimes_{\Bbb Q_p} \bold B_{\text{st}})),
\endalign
$$
where $\Bc$ and $\bold B_{\text{st}}$ are the rings of crystalline and semistable periods
(see \cite{F2}). Then $H^1_f(K,V) \simeq H^1_f(\Ddagrig (V))$ and $H^1_{st}(K,V)\simeq H^1_{st} (\Ddagrig (V)).$
\endproclaim
\demo{Proof} For any cocycle $x\in Z^1(G_K,V)$ denote by  $V_x$  the corresponding extension of $\Bbb Q_p$
by $V$. It is well known (and easy to check) that $\cl(x)\in H^1_f(K,V)$ if and only if
$\dim \Dc (V_x)=\dim \Dc (V)+1.$  Now the isomorphism $H^1_f(K,V) \simeq H^1_f(\Ddagrig (V))$ follows from the isomorphism
$$
\text{Ext}^1_{\Bbb Q_p[G_K]}(\Bbb Q_p(0),V) \simeq \text{Ext}^1_{\CR} (\Brigdag, \Ddagrig (V))
$$
given by Proposition 1.1.6 and Proposition 1.2.4. In the semistable case, the proof is analogous and is omitted here.
\enddemo

\proclaim{Lemma 1.4.3} Let $D$ be a $(\Ph,\Gamma)$-module. Then 
 $\cl (a,b)\in H^1_f(D)$ 
(resp. $\cl (a,b)\in H^1_{st}(D)$) 
if and only if
the equation $(\gamma-1)\,x=b$ has a solution in $D[1/t]$ 
(resp. in $D\otimes_{\Brigdag}\CR_{\log}(K) [1/t]$). 
\endproclaim
\demo{Proof}
An extension $D_\alpha$ is crystalline (resp. semistable) if and only if
there exists $x\in D[1/t]$ (resp. $x\in D\otimes_{\Brigdag}\CR_{\log}(K) [1/t]$)  such that 
$x+e \in \CDcris (D_\alpha)$ (resp. $x+e \in \Cal D_{\text{st}} (D_\alpha)$).  As $(\gamma -1)\,e=b,$ this proves the lemma.
\enddemo

The following proposition is proved (in a slightly different form) in \cite{FP}, Proposition 3.3.7 
and \cite{N}, sections 1.19-1.21. For the convenience of the reader
we recall the proof because it will be used in the proof of 
Proposition 1.5.8 below.

\proclaim{Proposition 1.4.4}  Let $D$ be a  potentially semistable $(\Ph,\Gamma)$-module. 
Then 

i) $H^0(D)$ and $H^1_f(D)$ are canonically isomorphic to the cohomology of the complex
$$
C^\bullet_{\text{\rm cris}} (D)\,\,:\,\, \CDcris (D)@>f>> t_D(K)\oplus \CDcris (D),
$$
where $t_D(K)=\CDdr (D)/\F^0 \CDdr (D)$ and $f(x)=(x\pmod {\F^0\CDdr (D)}, (1-\Ph)\,x).$

ii) $H^0(D)$ and $H^1_{st}(D)$ are canonically isomorphic to $H^i(C^\bullet_{\text{\rm st}} (D))$  ($i=0,1$) where
$$
C^\bullet_{\text{\rm st}} (D)\,\,:\,\,
\CDst (D)@>g>> t_D(K) \oplus \CDst (D)\oplus \CDst (D)@>h>> \CDst (D),
$$
with $g(x)=(x\pmod {\F^0\CDdr (D)}, (\Ph-1)\,x, N (x))$ and $h(x,y,z)= N(y)- (p\Ph -1)\,z).$
\endproclaim 
\demo{Proof} i)  Fix $r \gg 0$, $n \geqslant n(r)$ and consider the inclusion
$$
\iota_n\,\,:\,\, \CDcris (D)= (D^{(r)} [1/t])^{\Gamma} \hookrightarrow  (K_n((t))\otimes_{\bold B^{\dag, r}_{\text{rig},K},\,\iota_n} D^{(r)})^{\Gamma} \simeq \CDdr (D). 
$$
By Lemmas 5.1 and 5.4  of \cite{Ber1}, $\iota_n (x) \in \F^0 \CDdr (V)$ implies that $x\in D^{(r)}.$ Thus
$$
H^0(C^\bullet_{\text{\rm cris}} (D))= (\F^0\CDcris (D))^{\Ph=1} = (D^{(r)})^{\Ph=1,\gamma=1}= H^0(D).
$$
Next, let $D_{\alpha}$ be a crystalline extension of $D$. Then we have exact sequences
$$
\aligned
&0@>>>\CDcris (D) @>>>\CDcris (D_{\alpha}) @>>> K_0 @>>>0,\\
&0@>>>\Cal D_{\text{pst}} (D)@>>>\Cal D_{\text{pst}} (D_{\alpha}) @>>> K_0^{nr}@>>>0.
\endaligned
$$
Then $\CDcris (D_{\alpha})= \CDcris (D) \oplus K_0e$ where $b=(\Ph-1)\,e \in  \CDcris (D).$
The exact sequence
$$
0@>>>\F^0\CDdr (D)@>>>\F^0\CDdr (D_{\alpha}) @>>>K@>>>0
$$
shows that there exists $a \in \CDdr (D)$ such that $e+a \in \F^0\CDdr (D_{\alpha}).$ It is clear that
$a$ is unique modulo $\F^0\CDdr (D).$ If we replace $e$ by $e'=e+x$, $x\in \CDcris (D)$, then 
$b'=(\Ph-1)\,e'= b+ (\Ph-1)\,x$ and $a'=a-x.$ Thus the class of $(a\pmod {\F^0 \CDdr (D)}, b)$ modulo
$\text{Im} (f)$ does not depend on the choice of $d$ and gives a well defined element of 
$H^1 (C^\bullet_{\text{\rm cris}} (D)).$ 
Moreover  the filtered $(\Ph,N,G_K)$-module 
$\Cal D_{\text{pst}} (D_{\alpha})=\Cal D_{\text{pst}} (D)\oplus K_0^{nr} e$ is completely defined by
the class of $( a,b) \in t_D(K)\oplus \CDcris (D)$ modulo $\text{Im} (f)$ and the fact that $G_K$ acts trivially on $e$.
Conversely, to any   $(a\pmod {\F^0 \CDdr (D)}, b)\in t_D(K)\oplus \CDcris (D)$ we can associate the extension 
$M=\Cal D_{\text{pst}}(D)\oplus K_0^{nr}e$ of  filtered $(\Ph,N,G_K)$-modules defined by
$$\align 
&\Ph (e)=e+b, \quad N (e)=0,\\ 
&\F^i M_K=\cases \F^i \CDdr (D),&\text{\rm if $i < 0,$}\\
\F^i \CDdr (D) + K(e+a) , &\text{\rm if $i\geqslant 0.$}
\endcases
\endalign
$$
By Proposition 1.2.9 there exists a potentially semistable $(\Ph,\Gamma)$ module $D'$ such that
$\Cal D_{\text{pst}}(D')=M.$ Then $\CDcris (D')= M^{G_K, N=0}=\CDcris (D)\oplus K_0e .$ Thus
$D'$ is a crystalline extension of $\CR (K)$ by $D$  and i) is proved.

ii) In the semistable case the proof is analogous. First remark that
$$
H^0(C^\bullet_{\text{\rm st}} (D))= (\F^0\CDst (D))^{N=0,\Ph=1}=H^0(C^\bullet_{\text{\rm cris}} (D))=H^0(D).
$$
Next, if $D_{\alpha}$ is a semistable extension then we can write $\CDst(D_{\alpha})=\CDst (D)\oplus K_0e$
where $y=(\Ph-1)e\in \CDst (D)$ and $z=N(e)\in \CDst (D).$ As in i) there exists $x\in \CDdr (D)$ such that
$e+x\in \F^0\CDdr (D)$ and it is easy to see that the class of $(x \pmod{\CDdr (D)},y,z)$ modulo 
$\text{Im}(g)$ does not depend on the choice of $e$ and is a well defined element of $H^1(C^\bullet_{\text{\rm st}} (D)).$
Now ii) follows from Proposition 1.2.9.
\enddemo

\proclaim{Corollary 1.4.5} Let $D$ be a potentially semistable  $(\Ph,\Gamma)$-module. Then
$$
\align
&\dim_{\Bbb Q_p} H^1_f(D)-\dim_{\Bbb Q_p}H^0(D)\,=\,\dim_{\Bbb Q_p} t_D(K),\\
&\dim_{\Bbb Q_p} H^1_{st}(D)-\dim_{\Bbb Q_p}H^1_f(D)\,=\,\dim_{\Bbb Q_p} \CDcris(D^*(\chi))^{\Ph=1}.
\endalign
$$
\endproclaim
\demo{Proof} The first formula is obvious. The second follows from the fact that the cokernel 
of $h$ is dual to the kernel of $\CDst (D^*(\chi))@>(N,1-\Ph)>>\CDst (D^*(\chi))\oplus \CDst (D^*(\chi))$ which is
$\CDcris(D^*(\chi))^{\Ph=1}.$
\enddemo

\proclaim{Corollary 1.4.6} Let
$$
0@>>>D_1@>>>D@>>>D_2@>>>0
$$
be an exact sequence of potentially semistable $(\Ph,\Gamma)$-modules. Assume that one of the following
conditions holds:

a) $D$ is crystalline;

b) $\text{\rm Im} (H^0(D_2)@>>>H^1(D_1))\subset H^1_f(D_1).$

Then one has an exact sequence 
$$
0@>>>H^0(D_1)@>>>H^0(D)@>>>H^0(D_2)@>>>H^1_f(D_1)@>>>H^1_f(D)@>>>H^1_f(D_2)@>>>0 . \tag{3}
$$

\endproclaim
\demo{Proof} i) If $D$ is crystalline, then by Lemma 1.2.8  $D_1$ and $D_2$ are crystalline 
and we have an exact sequence of complexes
$$
0@>>> C^\bullet_{\text{\rm cris}} (D_1)@>>> C^\bullet_{\text{\rm cris}} (D) 
@>>> C^\bullet_{\text{\rm cris}} (D_2)@>>>0 .
$$
Passing to cohomology we obtain (3).  If the image of the connecting map is contained in $H^1_f(D_1)$
the sequence  (3) is again well defined. Only the exactness at $H^1_f(D_2)$ requires proof, but it
follows from the dimension argument using Corollary 1.4.5.  
\enddemo

\proclaim{Lemma 1.4.7} i) $D\mapsto H^i(C^\bullet_{\text{\rm st}} (D))$ defines a cohomological 
functor from $\bold{M}^{\Ph,\Gamma}_{{\text pst},K}$ to the category of $\Bbb Q_p$-vector spaces.
More precisely, let
$
0@>>>D_1@>>>D@>>>D_2@>>>0
$
be an exact sequence of $(\Ph,\Gamma)$-modules. 
If $D$ is potentially semistable, then $D_1$ and $D_2$ are potentially semistable and the sequence
$$
\multline
0@>>>H^0(D_1)@>>>H^0(D)@>>>H^0(D_2)@>>>H^1(C^\bullet_{st}(D_1))@>>>H^1(C^\bullet_{st}(D)) @>>>H^1(C^\bullet_{st}(D_2)) \\
@>>>H^2(C^\bullet_{st}(D_1))@>>>H^2(C^\bullet_{st}(D))@>>>H^2(C^\bullet_{st}(D_2))@>>>0
\endmultline 
$$
is exact. 

ii) The functor $D\mapsto H^i(C^\bullet_{\text{\rm st}} (D))$ is effacable. Namely for any
$\cl (x) \in H^i(C^\bullet_{\text{\rm st}} (D))$ ($i=1,2$) there exists an exact sequence
of semistable $(\Ph,\Gamma)$-modules $0@>>>D@>>>D'@>>>D''@>>>0$ such that
the image of $\cl (x)$ in $ H^i(C^\bullet_{\text{\rm st}} (D'))$ is zero.
\endproclaim
\demo{Proof} i)  By Lemma 1.2.8 $D_1$ and $D_2$ are potentially semistable. Let $L/K$ be a finite Galois extension such that
$D$, $D_1$ and $D_2$ are semistable over $L$ and let $I_{L/K}$ be the inertia subgroup of $\text{\rm Gal}(L/K).$
Then $H^1(I_{L/K}, \CDpst (D_1))=0$ and one has an exact sequence
$$
0@>>>\CDpst (D_1)^{I_{L/K}}@>>>\CDpst (D)^{I_{L/K}}@>>>\CDpst (D_2)^{I_{L/K}}@>>>0
$$
of  $K_0^{nr}$-vector spaces equipped with a semilinear action of $\text{\rm Gal} (K_0^{nr}/K_0).$
Taking invariants we obtain that the sequence 
$
0@>>>\CDst (D_1)@>>>\CDst (D)@>>>\CDst (D_2)@>>>0
$
is exact. Thus 
$$ 
0@>>> C^\bullet_{\text{\rm st}} (D_1)@>>> C^\bullet_{\text{\rm st}} (D) 
@>>> C^\bullet_{\text{\rm st}} (D_2)@>>>0 
$$
is an exact sequence of complexes. Passing to cohomology we obtain i).

ii) In dimension $1$, the effacability follows from the description of $H^1(C^\bullet_{\text{\rm st}} (D))$
in terms of extensions (see Proposition 1.4.4 ii)). Namely, let $\cl (x) \in H^1(C^\bullet_{\text{\rm st}} (D))$
and let $0@>>>D@>>>D_{x}@>>>\Brigdag @>>>0$ be the extension associated to $x$. Then the image of $\cl (x)$ in
$H^1(C^\bullet_{\text{\rm st}} (D_x))$ is zero.

Now prove that  $D\mapsto H^i(C^\bullet_{\text{\rm st}} (D))$ is effacable 
in dimension $2$. Let $x\in \CDst (D)$ and let $\cl (x)$ denote the class of $x$ in $H^2( C^\bullet_{\text{\rm st}} (D)).$
Set $x_1=Nx,$ $x_2=Nx_1, \ldots , x_m=Nx_{m-1}, Nx_m=0$ and consider the filtered $(\Ph,N, G_K)$-module $M'$ defined by
$$
\align
&M'=\CDpst (D) \oplus K_0e_1\oplus K_0e_2\oplus \cdots \oplus K_0 e_{m+1},\\
&Ne_i=e_{i+1}\quad\text{if $1\leqslant i\leqslant m$ and $Ne_{m+1}=0$},\\
&\Ph e_i=p^{-i}e_i +p^{-i}x_{i-1} \quad\text{if $1\leqslant i\leqslant m+1$},
\endalign
$$
$$
\F^i M_K'=\cases F^i\CDdr (D) \oplus Ke_1\oplus Ke_2\oplus \cdots \oplus K e_{m+1}
&\text{if $i\leqslant 0$}\\
F^i\CDdr (D) &\text{if $i>0$}.
\endcases  
$$
Then one has an exact sequence of filtered $(\Ph,N,G_K)$-modules
$$
0@>>>\CDpst (D)@>>>M'@>>>M''@>>>0 \tag{4}
$$
where $M''=\underset{i=1}\to{\overset{m+1}\to \oplus} K_0\overline e_i$, $N\overline e_i=\overline e_{i+1},$
$\Ph (\overline e_i)=p^{-i}\overline e_i$ and 
$$
\F^i M_K''=\cases M_K'' &\text{if $i\leqslant 0$}\\
0 &\text{if $i>0$}.
\endcases  
$$
By Proposition 1.2.9 the sequence (4) corresponds to an exact sequence of potentially
semistable $(\Ph,\Gamma)$-modules
$$
0@>>>D@>\alpha>>D'@>\beta>>D''@>>>0
$$
such that $M'=\CDpst (D')$ and $M''=\CDpst (D'')$. 
Since $(p\,\Ph-1)\,e_1=x,$ the image of $\cl (x)$ in $H^2 (C^\bullet_{\text{\rm st}} (D'))$ is zero and the lemma is proved.
\enddemo

{\bf 1.4.8.} If $D_1$ and $D_2$ are two semistable modules we define a pairing
$$
\cup \,\,:\,\,H^i(C^\bullet_{\text{\rm st}} (D_1))\times  H^{j}(C^\bullet_{\text{\rm st}} (D_2)) @>>>H^{i+j} (C^\bullet_{\text{\rm st}}(D_1\otimes D_2))
$$
by the formulas
$$
\cl (x_1) \cup \cl (x_2)=\cl (x_1\otimes x_2) \qquad \text{if $i=0$, $j=0$},
$$
$$
\multline
\cl (x_1) \cup \cl (x_2\pmod{\F^0\CDdr (D_2)},y_2,z_2)\,=\,\cl (x_1\otimes  x_2\pmod{ \F^0{\CDdr (D_1\otimes D_2)}},x_1\otimes y_2,x_1\otimes z_2),\\
\text{if $i=0$, $j=1$},
\endmultline
$$
$$
\cl (x_1) \cup \cl (x_2)=\cl (x_1\otimes x_2) \qquad \text{if $i=0$, $j=2$},
$$
$$
\multline \cl (x_1\pmod{\F^0\CDdr (D_2)},y_1,z_1) \cup \cl (x_2\pmod{\F^0\CDdr (D_2)},y_2,z_2)=\cl (z_1\otimes y_2 +y_1\otimes p\Ph (z_2))\\
\qquad \text{if $i=1$, $j=1$}.
\endmultline
$$
An easy computation shows that this product is compatible with connecting homomorphisms. Namely,
if $0@>>>D_1@>>>D_1'@>>>D_1''@>>>0$ is exact, then the diagrams
$$
\CD
H^i(C^\bullet_{\text{\rm st}} (D_1''))\times H^j(C^\bullet_{\text{\rm st}} (D_2))@>\cup >>H^{i+j}(C^\bullet_{\text{\rm st}} (D_1''\otimes D_2))\\
@VVV @VVV \\
H^{i+1}(C^\bullet_{\text{\rm st}} (D_1))\times H^j(C^\bullet_{\text{\rm st}} (D_2))@>\cup >>H^{i+j+1}(C^\bullet_{\text{\rm st}} (D_1\otimes D_2))
\endCD
$$
commute.
\proclaim{Proposition 1.4.9} There exists a unique natural transformation of cohomological functors 
\linebreak
$h^*\,:\,H^*(C^\bullet_{\text{\rm st}} (D)) @>>> H^*(D)$
satisfying the following properties:

1) $h^0$ and $h^1$ coincide with the maps $H^0(C^\bullet_{\text{\rm st}} (D))\simeq H^0(D)$  and 
$H^1(C^\bullet_{\text{\rm st}} (D))\simeq H^1_{st}(D)@>>>H^1(D)$  given by Proposition 1.4.4.

2) $h^*$ is compatible with cup-products. 
\endproclaim
\demo{Proof} Remark that $h^0$ and $h^1$ are already defined. Let $\cl (x) \in  H^2(C^\bullet_{\text{\rm st}} (D))$.
We define $h^2(\cl (x))$ by the usual way using  effacability. Let 
$$
E_D\,\,:\,\,0@>\alpha>>D@>>>D'@>\beta>>D''@>>>0 
$$
be an exact sequence in ${\bold M}^{\Ph,\Gamma}_{\text{pst},K}$ such that $\cl (x)$ vanishes in 
$H^2(C^\bullet_{\text{\rm st}} (D))$.
Consider the diagram
$$
\xymatrix{
\cdots\ar[r] &H^1(D') \ar[r] &H^1(D'') \ar[r]^{\delta_1}  & H^2(D) \ar[r]   &H^2(D')\ar[r] & \cdots\\
\cdots\ar[r] &H^1(C^\bullet_{\text{\rm st}} (D')) \ar[r] \ar[u] &H^1(C^\bullet_{\text{\rm st}} (D''))\ar[r]^{\Delta_1} \ar[u]^{h^1}
&H^2(C^\bullet_{\text{\rm st}} (D))\ar[r]    &H^2(C^\bullet_{\text{\rm st}} (D'))\ar[r]&\cdots
}
$$
Since  the image of $\cl (x)$ in $H^2 (C^\bullet_{\text{\rm st}} (D'))$ is zero, there exists $a\in 
H^1 (C^\bullet_{\text{\rm st}} (D''))$ such that $\Delta_1(a)=\cl (x).$
Set  $h^2(\cl (x))=\delta_1\circ h^1 (a).$ It is well known (see for example \cite{Sz}) that for abelian
categories this  gives a well defined morphism $H^2(C^\bullet_{\text{\rm st}} (D))@>>> H^2(D).$
In our case this can be proved by the same way. Namely, for any morphism
$f\,\,:\,\,D@>>>P$ of potentially semistable $(\Ph,\Gamma)$-modules the usual construction 
gives an exact sequence $f\circ E_D$ which seats in the diagram:
$$
\xymatrix{
E_D: &0 \ar[r] &D\ar[r]^{\alpha} \ar[d]^f &D' \ar[r]^{\beta} \ar[d] &D'' \ar[r] \ar[d]^{=} &0 \\
f\circ E_D: &0\ar[r] &P\ar[r] &P' \ar[r] &D'' \ar[r] &0.}
$$ 
Indeed, here $P'=\text{coker} (D@>(\alpha ,f)>> D'\oplus P)$ is a free $\Brigdag$-module
because $D$ is saturated in $D'\oplus P$. Let $E'_D\,\,:0@>>>D@>>>X'@>>>X''@>>>0$ be another
exact sequence such that $\cl (x)$ vanishes in $H^2(C^\bullet_{\text{\rm st}} (X')).$
Set $f_2\,\,:\,\,D\oplus D@>>>P$, $f_2(d_1,d_2)=d_1+d_2$ and consider $f_2\circ (E_D\oplus E_D').$
The  injections $i_{1,2}\,\,:\,\,D @>>>D\oplus D$, $i_1(d)=(d,0)$, $i_2(d)=(0,d)$   
induce  natural morphisms $E_D @>>> f_2\circ (E_D\oplus E_D')$ and 
$E'_D @>>> f_2\circ (E_D\oplus E_D').$ This allows to show that $h_2(\cl (x))$ is well defined. A similar
argument proves that $h_2$ is a morphism. The uniqueness of $h_2$ is clear. It remains to show
that the diagrams
$$
\CD
H^i(C^\bullet_{\text{\rm st}} (D_1))\times H^j(C^\bullet_{\text{\rm st}} (D_2))@>\cup >>H^{i+j}(C^\bullet_{\text{\rm st}} (D_1\otimes D_2))\\
@VVV @VVV \\
H^i(D_1)\times H^j(D_2)@>\cup>>H^{i+j}(D_1\otimes D_2)
\endCD 
$$
commute. It is immediate   if  $i=j=0$  and the general case can be deduced using dimension shifting.

\enddemo

\proclaim{Corollary 1.4.10} Let $D$ be a potentially semistable  $(\Ph,\Gamma)$-module. Then $H^1_f(D^*(\chi))$
is the orthogonal complement to $H^1_f(D)$ under the duality
$$
H^1(D)\times H^1(D^*(\chi)) @>\cup >> H^2 (\Brigdag (\chi))\simeq \Bbb Q_p .
$$ 
\endproclaim 
\demo{Proof} By Proposition 1.4.9 we have a commutative diagram
$$
\CD
H^1(C^\bullet_{\text{\rm st}}(D))\times H^1(C^\bullet_{\text{\rm st}}(D^*(\chi))) @>\cup >> H^2(C^\bullet_{\text{\rm st}}(\Brigdag (\chi)))\\
@VVV @VVV\\
H^1(D)\times H^1(D^*(\chi)) @>\cup>>\Bbb Q_p
\endCD
$$
From the definition of the cup product it follows immediately that $x\cup y=0$ for all $x\in H^1(C^\bullet_{\text{\rm cris}}(D))$ and $y\in H^1(C^\bullet_{\text{\rm cris}}(D^*(\chi)))$. This proves that $H^1_f(D^*(\chi))$ and $H^1_f(D)$ are orthogonal to each other.
Next,  by Corollary 1.4.5 together with the Euler-Poincar\'e characteristic formula we have
$$
\dim_{\Bbb Q_p}H_f^1(D)+\dim_{\Bbb Q_p} H^1_f(D^*(\chi))=\dim \CDdr (D)+\dim_{\Bbb Q_p}H^0(D) +\dim_{\Bbb Q_p}H^0(D^*(\chi))=\dim_{\Bbb Q_p}H^1(D).
$$
The corollary is proved.
\enddemo

{\bf 1.5. Semistable modules of rank 1.}

{\bf 1.5.1.} In this section we compute  $H^1_f(D)$ of semistable modules of rank $1$.  For any continuous
character $\delta \,:\,\Bbb Q_p^*@>>> \Bbb Q_p^*$ let $\CR (K) (\delta)$
denote the $(\Ph,\Gamma)$-module $\Cal R (K) e_{\delta}$ such that
$\Ph (e_{\delta})=\delta (p) e_{\delta}$ and 
$\gamma (e_{\delta})=\delta (\chi(\gamma))\,e_{\delta},$ $\gamma \in \Gamma$. Write $ x$ for
the character given by the identity map
and $\vert x\vert$ for $\vert
x\vert=p^{-v_p(x)}.$

\proclaim{Lemma 1.5.2} The following statements are equivalent:

i) $D$ is a semistable module of rank $1$;

ii) $D$ is a crystalline module of rank $1$;

iii) $D$ is isomorphic to $\CR(K) (\delta)$ where $\delta \,:\,\Bbb Q_p^* @>>>\Bbb Q_p^*$ is such that
$\delta (u)=u^m$ ($u\in \Bbb Z_p^*$) for some $m\in \Bbb Z$.
\endproclaim
\demo{Proof}   Since the operator $N\,:\, \CDst(\CR (K)(\delta))@>>>\CDst(\CR(K)(\delta))$ is nilpotent,  it is clear that $i) \Leftrightarrow ii)$ and $iii)\Rightarrow i).$ We prove that $i) \Rightarrow iii).$  If  $\CR (K)(\delta)$ is semistable
it is  a de Rham and there exists 
$n\in \Bbb Z$ such that $\delta (u)=u^m$ on  $1+p^n\Bbb Z_p\subset \Bbb Z_p^*$. Replacing $\delta$ by $\delta x^{-m}$ we may assume that
$m=0.$ As $(\CR (K)[1/t])^{\Gamma_n} =K_0,$ we obtain that $\CDcris (\CR (K)(\delta)) =(K_0(\delta))^{\Gamma}$ is a one-dimensional $L$-vector space.
Thus $\delta=0$ and the proposition is proved.
\enddemo

\proclaim{Proposition 1.5.3} Let   $\CR(K)(\delta)$ be such that $\delta (u)=u^m$ ($u\in \Bbb Z_p^*)$ for some $m\in \Bbb Z$. Then

i) $$H^0(\CR (K)(\delta))= \cases  \Bbb Q_p t^m &\text{\rm if
$\delta=x^{-m}$, $k\in \Bbb N$}\\
0 & \text{\rm otherwise.}
\endcases
$$

ii) Assume that  $m\geqslant 1$. If $\delta \ne \vert x\vert x^m$ then $H^1 (\CR (K)(\delta))=H^1_f(\CR (K)(\delta))$
is a $\Bbb Q_p$-vector space of dimension $[K:\Bbb Q_p].$ If $\delta = \vert x\vert x^m$ then 
$H^1_{st}(\CR (K)(\delta))=H^1(\CR (K) (\delta))$ is a $\Bbb Q_p$-vector space of dimension
$[K:\Bbb Q_p]+1$ and $\dim_{\Bbb Q_p}H^1_f(\CR (K) (\delta))=[K:\Bbb Q_p].$

iii) Assume that $m\leqslant 0.$ If $\delta \ne x^m$, then $H^1(\CR (K)(\delta))$ is a $\Bbb Q_p$-vector space
of dimension $[K:\Bbb Q_p]$ and $H^1_f(\CR (K)(\delta))=0.$ If $\delta=x^m$, then 
$\dim_{\Bbb Q_p}H^1(\CR (K) (\delta))=[K:\Bbb Q_p]+1$ and $\dim_{\Bbb Q_p}H^1_f(\CR(K) (\delta))$ is the
one dimensional $\Bbb Q_p$-vector space generated by $\cl (t^{-m},0)\otimes e_{\delta}.$
\endproclaim
\demo{Proof} This is a direct application of the Euler-Poincar\'e caracteristic formula
together with Corollary 1.4.5.
\enddemo

{\bf 1.5.4.}  In the remainder of this paragraph we suppose that $K=\Bbb Q_p$ and
write $\CR$ for $\CR (\Bbb Q_p)$. 
In \cite{Cz4}, Colmez  studies  general  $(\Ph,\Gamma)$-modules of rank $1$ over 
$\CR$. He proves that any $(\Ph,\Gamma)$-module of rank $1$
is isomorphic to $\CR(\delta)$ for some $\delta \,:\,\Bbb Q_p^*@>>>\Bbb Q_p^* $ (see \cite{Cz4}, Proposition 3.1)
and describes $H^1 (\CR(\delta))$ in details.

\proclaim{Proposition 1.5.5} i)  If $\delta  =x^{-m},$ $m\geqslant 0$ then   
 $\cl (t^m,0)\,e_{\delta} $ and $\cl (0,t^m)\,e_{\delta}$ form a $L$-basis of 
$H^1(\CR(\delta))$.

iv) If $\delta = \vert x\vert x^{m},$ $m\geqslant 1,$ then  
$H^1(\CR (\delta))$ is generated over $\Bbb Q_p$ by $\cl (\alpha_m)$ and $\cl (\beta_m)$
where
$$
\align
&\alpha_m= \frac{(-1)^{m-1}}{(m-1)!}\, \partial^{m-1} \left (\frac{1}{\pi}+\frac{1}{2},a \right )\, e_{\delta} ,\qquad
(1-\Ph)\,a=(1-\chi (\gamma) \gamma)\,\left (\frac{1}{\pi}+\frac{1}{2} \right ),\\
&\beta_m=\frac{(-1)^{m-1}}{(m-1)!}\,
\partial^{m-1} \left (b,\frac{1}{\pi} \right )\,e_{\delta}, \qquad
(1-\Ph)\,\left (\frac{1}{\pi}\right )\,=\,(1-\chi (\gamma)\,\gamma)\,b
\endalign
$$
and $\partial $ denotes the differential operator $(1+\pi)\,d/d\pi .$
\endproclaim
\demo{Proof} See sections 2.3-2.5 of \cite{Cz4}.
\enddemo

To simplify notations we will write $e_m$ (respectively $w_m$) for $e_\delta$
if $\delta= \vert x\vert x^m$, $m\geqslant 1$ (respectively if $\delta =x^{-m}$, $m\geqslant 0$).

\proclaim{Corollary 1.5.6}
Any element $\cl ((f,g)e_m)\in H^1(\CR (x^m\vert x\vert))$ can be written in the form
$$
\cl \left((f,g)e_m \right )\,=\,x \, \cl (\alpha_m)\,+\,y \, \cl (\beta_m)
$$
where $x= \res (ft^{m-1}dt),$ $y=\res (gt^{m-1}dt).$

\endproclaim
\demo{Proof} The case $m=1$ is implicitly contained in the proof of Proposition 2.8 of \cite{Cz4}
and the general case can be proved by the same argument. Namely, the formulas
$$
\res (\chi (\gamma)\gamma (h) dt)\,=\,\res (hdt),\qquad
\res (\Ph (h)dt)\,=\,\res (hdt)
$$
imply that 
$$
\res (\chi (\gamma)^m \gamma (h)\,t^{m-1}dt)\,=\,
\res (ht^{m-1}dt),\qquad \res (p^{m-1}\Ph (h) t^{m-1} dt)\,=\,\res (ht^{m-1}dt).
$$
Thus  $ \res (ft^{m-1}dt)$ and  $\res (gt^{m-1}dt)$ do not change 
if we add to $(f,g)e_m$ a boundary 
\linebreak
$((1-\Ph)\,(he_m),(1-\gamma)\,(he_m)).$ As $a\in \Cal R^+=\CR \cap \Bbb Q_p[[\pi]]$ (see \cite{Ben1}, Lemma 2.1.3 or \cite{Cz4}),
we have $\res (\alpha_m t^{m-1} dt)\,=\,(1,0).$  For any $g\in \A_{\Bbb Q_p}^+$ set
$$
\ell (g)\,=\,\frac{1}{p}\log \left (\frac{g^p}{\Ph (g)}\right ).
$$
As $g^p/\Ph (g)\equiv 1 \pmod{p}$ in $\bold A_{\Bbb Q_p}$, it is easy to see that
$\ell (g)\in \bold A_{\Bbb Q_p}.$ Let $\psi$ denote the operator 
$\psi (G(\pi))=\dsize\frac{1}{p}\Ph^{-1}\underset{\zeta^p=1}\to \sum G((1+\pi)\zeta -1).$ 
A short computation shows that
$\psi (\ell (\pi))=0.$ Then by Lemma 1.5.1 of \cite{CC} there exists a unique  $b_0\in \bold A_{\Bbb Q_p}^\dagger$ such that 
$(\gamma -1)\,b_0\,=\,\ell (\pi).$ Applying the operator $\partial$ to this formula we obtain that
$b=\partial b_0.$ Thus
$\res ((\partial^{m-1}b)t^{m-1}dt)\,=\,0$ and $\res (\beta_mt^{m-1}dt)=(0,1).$
By Proposition 1.5.5 any  $\cl ((f,g)w_m)\in H^1(\CR_L(\vert x\vert x^m)$ can be written in the form 
$
\cl \left((f,g)w_m \right )\,=\,x \, \cl (\alpha_m)\,+\,y \, \cl (\beta_m)
$
Applying previous formulas we obtain  that $x= \res (ft^{m-1}dt)$ 
and $y=\res (gt^{m-1}dt).$
The corollary is proved.
\enddemo

{\bf 1.5.7.} Set 
$$
\align
&\alpha_m^* =\left (1-\frac{1}{p}\right )\,\cl (\alpha_m),\quad \beta_m^*=\left (1-\frac{1}{p}\right )\log \chi (\gamma) \,\cl (\beta_m)
\qquad \text{if  $m \geqslant 1$},\\
&x_m^*=\cl (t^m,0) w_m, \quad y_m^*=\log \chi (\gamma)\,\cl (0,t^m) w_m \qquad \text{if  $m\geqslant 0$.}
\endalign
$$
If $m=0$, then $H^1(\CR_L (x^{m})=H^1(\CR))$ is canonically isomorphic to $H^1(\Bbb Q_p,\Bbb Q_p)=
\Hom (G_{\Bbb Q_p},\Bbb Q_p)$. Let $\text{ord} \,:\,\text{\rm Gal} (\Bbb Q_p^{ur}/\Bbb Q_p)@>>>\Bbb Z_p$ denote the character
determined by $\text{ord} (\Ph)=1$. Then $\text{ord}$ and $\log \chi$ form a canonical basis of 
$\Hom (G_{\Bbb Q_p},\Bbb Q_p)$ which corresponds to  $x_0^*$, $y_0^*$. If $m=1$, then $H^1(\CR (\vert x\vert x))\simeq H^1(\CR (\chi))$
is isomorphic to $H^1(\Bbb Q_p,\Bbb Q_p(1)).$ Let $\kappa \,:\,\Bbb Q_p^*@>>>H^1(\Bbb Q_p,\Bbb Q_p(1))$ denote the Kummer map.
Then $\kappa (u)=-\log (u)\alpha_1^*$  if $u\equiv 1 \pmod p$ and $\kappa (p)=- \beta_1^*$ (see \cite{Ben1}, Proposition 2.1.5).

\proclaim{Proposition 1.5.8} i) Let $m\geqslant 1$. The basis $\alpha_m^*$, $\beta_m^*$
is dual to the basis  $y_m^*$,  $x_m^*$ under the pairing
$$
H^1(\CR (\vert x\vert x^m))\times H^1(\CR (x^{1-m}))@>\cup>> \Bbb Q_p .
$$

In particular, $\alpha_m^*$ generates  $H^1_f(\CR (\vert x\vert x^m)).$

ii) The canonical isomorphism
$$
H^1(\CR (\vert x\vert x^m)) \iso H^1 (C^\bullet_{\text {\rm st}}(\CR (\vert x\vert x^m)))
$$
sends $\alpha^*_m$ to $-\cl ( 1,0,0)$ and $\beta^*_m$ to $-\cl (0,0,1).$

\endproclaim 
\demo{Proof} i) The proof follows from the construction of the cup product (see 1.1.5) together
with the explicit description of the isomorphism $H^2(\CR (\chi))\simeq \Bbb Q_p$ reviewed in 1.1.7.
Remark that 
$$
\partial^{m-1} \left (\frac{1}{\pi}+\frac{1}{2}\right ) \equiv \frac{(-1)^{m-1}(m-1)!}{t^m}
\pmod{\Bbb Q_p[[\pi]]} .
$$
Thus
$$
\alpha_m^* \cup y_m^* = \frac{(-1)^{m-1}}{(m-1)!}\res \left (\partial^{m-1}\left (\frac{1}{\pi}+\frac{1}{2} \right ) t^{m-1}\,\frac{d\pi}{1+\pi}\right ) =1.
$$
On the other hand it is clear that $\alpha_m^*\cup x_m^*=0.$
The proof of other formulas is analogous.

iia) Recall that $e_m$ denote the canonical basis of $\CR (\vert x\vert x^m).$
Let 
$$
0@>>>\CR (\vert x\vert x^m) @>>>D@>>>\Cal R@>>>0
$$
be the extension associated to $\alpha_m .$ Then $D=\CR(\vert x\vert x^m) +\Cal R u_m$ where
$$
(\Ph-1)\,u_m\,=\,\frac{(-1)^{m-1}}{(m-1)!}\,\,\partial^{m-1} \left (\frac{1}{\pi}+\frac{1}{2} \right ) e_m,\qquad (\gamma -1)\,u_m\,=\,\frac{(-1)^{m-1}}{(m-1)!}\,\, \partial^{m-1}(a)\, e_m.
$$
As $\cl (\alpha_m)\in H^1_f(\CR (\vert x\vert x^m),$ we have an exact sequence
$$
0@>>>\Cal D_{\text{\rm cris}}(\CR (\vert x\vert x^m))@>>> \Cal D_{\text{\rm cris}}(D)@>>>
\Cal D_{\text{\rm cris}}(\Cal R)@>>>0 \tag{5}
$$
which induces an isomorphism
$$
\Cal D_{\text{\rm cris}}(D)^{\Ph=1}\simeq \Cal D_{\text{\rm cris}}(\Cal R)^{\Ph=1}.
$$
Then there exists a unique $c_me_m\in \CR (\vert x\vert x^m)\left [1/t\right ]$ such that
$$
(\Ph-1)\,(u_m+ c_me_m)\,=\,0.
$$
Thus
$
(1-p^{m-1}\Ph )\,c_m\,=\,-\dsize\frac{(-1)^{m-1}}{(m-1)!}\,\,\partial^{m-1} \left (\frac{1}{\pi}+\frac{1}{2} \right )
$
and a short computation shows that 
\linebreak
$c_m=\dsize\frac{(-1)^{m-1}}{(m-1)!}\,\,\partial^{m-1}c_1.$
For $m=1$ this equation can be written in the form
$$(1-\Ph/p)\,(tc_1)=-t\left (\dsize\frac{1}{\pi}+\frac{1}{2} \right )\in \Cal R^+,$$
where $\CR^+ =\CR \cap \Bbb Q_p[[\pi]]$ and  $t\left (\dsize\frac{1}{\pi}+\frac{1}{2} \right )\,\equiv \,1 \pmod{\pi^2}.$
Then Lemma A.1 of \cite{Cz4} implies that $tc_1\in \Cal R^+$ and therefore satisfies
$tc_1 \equiv -\left (1-\frac{1}{p} \right )^{-1} \pmod{\pi}.$ Taking derivations 
we obtain by induction that $t^m c_m \equiv -\left (1-\frac{1}{p} \right )^{-1} \pmod{\pi}.$ 
Now from the proof of Proposition 1.4.4 it follows that the extension (5) corresponds to the class
$-\left (1-\frac{1}{p} \right )^{-1}  \cl (1,0,0) \in H^1(C_{\text{st}}^\bullet (\CR (\vert x\vert x^m)))$
and the first formula of ii) is proved. 


iib) Let $0@>>>\CR (\vert x\vert x^m)@>>>D@>>>\CR @>>>0 $ be the extension
associated to $\cl (\beta_m)$. Write $D=\CR (\vert x\vert x^m) +\CR v_m$ where
$$
(\Ph-1)\,v_m\,=\,\frac{(-1)^{m-1}}{(m-1)!}\,\,\partial^{m-1} (b)\,e_m,\qquad 
(\gamma-1)\,v_m\,=\,\frac{(-1)^{m-1}}{(m-1)!}\,\,\partial^{m-1}\left (\frac{1}{\pi}\right )\,e_m.
$$
By Proposition 1.5.3 $D$ is semistable and 
there exists $d_me_m\in \CR_{\log} (\vert x\vert x^m) [1/t]$
such that 
$(\gamma-1)\, (v_m+d_me_m)\,=\,0.$ Thus
$$
(\chi (\gamma)^m \gamma -1)\,d_m\,=\,\frac{(-1)^{m}}{(m-1)!}\,\,\partial^{m-1}\left (\frac{1}{\pi}\right ).
$$
An easy computation shows that $d_m=\dsize\frac{(-1)^{m-1}}{(m-1)!}\,\,\partial^{m-1}\left (d_1 \right ).$
Set $\nabla =\log (\gamma)/\log \chi (\gamma)$ and 
\linebreak
$\nabla_0=\nabla /(\gamma-1).$
As  $\nabla= t\partial $  on $\CR$ 
we see that one can take $d_1=-t^{-1}\nabla_0 (\log \pi)+(\chi (\gamma)-1)^{-1}.$ From
$$
(\Ph-1)\,d_1=t^{-1} \nabla_0 \left (1-\frac{\Ph}{p} \right )\log \pi,
$$
it follows that
$$
(\gamma-1)\,\left ((\Ph-1)\,d_1e_1 \right )\,=\,\nabla \left (1-\frac{\Ph}{p} \right )( \log \pi) e_1= (1-\Ph) \,\partial (\log \pi)e_1=
(1-\Ph) \left (\frac{1}{\pi}\right ) e_1.
$$
As $\psi (tb)=\dsize\frac{1}{p}\psi (\Ph (t)b)=t\psi (b)/p=0,$ the element $bt\in \CR^{\psi=0}$ is a solution
of the equation $(\gamma-1)\,x=(1-\Ph) (1/\pi).$ On the other hand, as
$
\underset {\zeta^p=1}\to \prod ((1+X)\,\zeta -1)=\Ph (X),
$
we have 
$$
\psi \log \left (\frac{\Ph (X)}{X^p}\right )\,=\,\frac{1}{p}\, \Ph^{-1} \log \left (
\frac{\Ph (X)^p}{\underset {\zeta^p=1}\to \prod ((1+X)\,\zeta -1)^p}\right )\,=\,0.
$$ 
Thus $t(\Ph-1) d_1=(1-\Ph)\nabla_0 \left (1-1/p \right )\,\log \pi \in \CR^{\psi=0}$ is also
a solution of $(\gamma-1)\,x=(1-\Ph) (1/\pi).$ As $\gamma-1$ is bijective on $\CR^{\psi=0}$ 
(see for example \cite{Ber3}, Lemma I.3) this implies that
 $(\Ph-1)\,d_1=-b$ and $(\Ph-1)\,(d_me_m)=\frac{(-1)^{m}}{(m-1)!}\,\,\partial^{m-1}(b)\, e_m.$
Using the formula  $\partial\,\gamma=\chi (\gamma)\,\gamma \,\partial$ we obtain 
$$
\partial^{m-1} \left (t^{-1}\nabla_0 (\log \pi)\right )\,=\,(-1)^{m-1} (m-1)! \, t^{-m}\nabla_0(\log \pi)+t^{1-m} x,
$$
where $x\in \pi^{1-m} \CR^+$.
Set $\gamma_1=\gamma^{p-1}$. Then
$$
\nabla_0 =\left ( \sum_{i=0}^{p-2} \gamma^i \right )\,\frac{\nabla}{\gamma_1-1}\,=\,(\log \chi (\gamma_1))^{-1} 
\left (\sum_{i=0}^{p-2} \gamma^i \right )\,
 \underset{n=0}\to {\overset \infty\to \sum} (-1)^n \dsize\frac{(\gamma_1-1)^n}{n+1}. 
$$
Remark that
$$
\left ( \sum_{i=0}^{p-2} \gamma^i \right )\log \pi=(p-1) \log \pi+ y_1
$$
where 
$$
y_1=\sum_{i=1}^{p-2} (\gamma^i-1)\log \pi =\log \left (\frac{\pi^{1+\gamma+\cdots +\gamma^{p-2}}}{\pi^{p-1}}\right ).
$$
In particular, 
$$
\iota_1 (y_1) (0)=\log \left (\frac{N_{\Bbb Q_p(\zeta_p)/\Bbb Q_p} (\zeta_p-1)}{(\zeta_p-1)^{p-1}}\right )\,=\,
\log \left (\frac{p}{(\zeta_p-1)^{p-1}}\right )\,=\,-(p-1)\log (\zeta_p-1).
$$
Set 
$$
y_2=\left ( \sum_{i=0}^{p-2} \gamma^i \right ) \underset{n=1}\to {\overset \infty\to \sum} (-1)^n \dsize\frac{(\gamma_1-1)^n}{n+1}
\log \pi.
$$
As $\iota_1 ((\gamma_1-1)\log \pi) (0)=0$ we have $\iota_1 (y_2) (0)=0$ and
$$
d_m= - \,(\log \chi (\gamma))^{-1}\,t^{-m} (\log \pi + y)
$$
where $y\in \bold B^{\dagger,p-1}_{\text{\rm rig},\Bbb Q_p}$ is such that $\iota_1 (y)=-\log (\zeta_p-1).$
Remark that as usually the $p$-adic logarithm is normalized by $\log (p)=0.$
Thus
$$
\cl(\beta_m)\,=\,(\log \chi (\gamma))^{-1}\left (1-\frac{1}{p}\right ) t^{-m}\cl \left ((p^{-1}\Ph-1) (\log \pi +y), (\gamma-1)(\log \pi+y) \right ) e_m .
\tag{6}
$$ 
Let $G(\vert x\vert,m)$ and $G'(\vert x\vert,m)$ be two elements of $\bold B^{\dagger,p-1}_{\text{\rm rig},\Bbb Q_p}$
satisfying the following conditions:
$$
\iota_n (G(\vert x\vert ,m)) \equiv p^{-n} \pmod {t^m},
\quad
\iota_n (G'(\vert x\vert ,m)) \equiv \log (\zeta_{p^n}e^{t/p^n}-1) \pmod {t^m},
\quad  \forall n\geqslant 1.
$$

By \cite{Cz4}, Proposition 2.19 there exist unique $\lambda,\mu \in \Bbb Q_p$ such that
$\beta^*_m=\cl (a,b) e_m$ where 
$$
\aligned
&a=t^{-m}(p^{-1}\Ph-1) (\lambda G(\vert x\vert,m)+\mu (\log \pi -G'(\vert x\vert,m))),\\
&b=t^{-m}(\gamma-1) (\lambda G(\vert x\vert,m)+\mu (\log \pi -G'(\vert x\vert,m))).
\endaligned
\tag{7}
$$     
Remark that if two cocycles  $t^{-m}((p^{-1}\Ph-1) x, (\gamma-1) x)$ and    
 $t^{-m}((p^{-1}\Ph-1) y, (\gamma-1) y)$ ($x,y\in \CR_{\log})$ are homologous
 in $Z^1(\CR (\vert x\vert x^m),$ then  $y=x+t^{m}z$ for some $z\in \CR .$ In particular,
 if $x,y\in \bold B^{\dagger,p-1}_{\text{\rm rig},\Bbb Q_p}$ then $\iota_1(x) (0)=\iota_1 (y) (0).$
Comparing (6) and (7) we find that $\lambda =0$ and 
$\mu= (\log \chi (\gamma))^{-1}.$ Therefore there exists a lifting $v_m'\in D$ of $1\in \CR$ such that 
the $(\Ph,N)$-module $\CDst (D)=(D\otimes_{\CR}\CR_{\log}[1/t])^{\Gamma}$ is generated by
$f_1= t^{-m} e_m$ and $f_2= v_m'+(\log \chi (\gamma))^{-1}\,(\log \pi -G'(\vert x \vert ,m)) f_1.$ We see 
immediately that $\Ph (f_2)= f_2$ and $N(f_2)=-(1-1/p)^{-1}(\log \chi (\gamma))^{-1}f_1$ and $f_2\in \F^0 \CDst (D).$
Thus $D$ corresponds to  the class $(1-1/p)^{-1}(\log \chi (\gamma))^{-1}\cl (0,0,-1) \in H^1(C^\bullet_{\text st} (\CR (\vert x\vert x^m)).$
The proposition is proved.
\enddemo

\proclaim{Proposition 1.5.9} Let $D$ be a semistable $(\Ph,\Gamma)$-module of rank $d$ 
with Hodge-Tate weights  $k_1,\ldots, k_d$. Assume that $\CDst (D)=\CDst (D)^{\Ph=\lambda}$ 
for some $\lambda \in \Bbb Q_p.$
 Then 
$$
D \simeq \underset{i=1}\to{\overset d\to \oplus}\CR(\delta_i),
$$
where $\delta_i$ are defined by $\delta_i(u)=u^{-k_i}$ ($u\in \Bbb Z_p^*$) and $\delta_i(p)=\lambda p^{-k_i}.$ 
\endproclaim
\demo{Proof} 1) We prove the Proposition by induction on $d=\text{rg} (D).$ The case $d=1$ is trivial. 
Let $D$ be a semistable $(\Ph,\Gamma)$-module of rank $2$ with Hodge weights $k_1\leqslant k_2.$ Choose
a non-zero $v\in \F^{k_2}\CDst (D)$ and put $\Cal F_1\CDst (D)=\Bbb Q_p v$, $\Cal F_2\CDst (D)=D.$ By Proposition 1.3.2 
the triangulation of $D$ associated to this filtration gives rise to an exact sequence
$$
0@>>>\CR(\delta_2)@>>>D@>>>\CR(\delta_1)@>>>0
$$
where $\delta_i(p)=\lambda p^{-k_i}$ and $\delta_i(u)=u^{-k_i}$, $u\in \Bbb Z_p^*.$  Put $\delta=\delta_2\delta_1^{-1}.$
Then $\delta (x)=x^{-k}$ with $k=k_2-k_1\geqslant 0$ and $D(\delta_1^{-1})$ can be viewed as a semistable extension
$$
0@>>>\CR (x^{-k})@>>>D(\delta_1^{-1}) @>>>\CR @>>>0.
$$  
The relation $N\,\Ph=p\,\Ph N$ together with the fact that $\CDst (D)$ is pure implies that 
$\CDst (D)^{N=0}=\CDst (D).$ Thus $D$ is crystalline and by Proposition 1.5.3 the class of $D(\delta^{-1})$ in
$H^1(\CR (\delta))$ is $ ax_k^*$ for some $a\in \Bbb Q_p.$ Write $D(\delta_1^{-1})=\CR e_{\delta}\oplus \CR e$
where $e\in D$ is the lifting of $1\in \CR $ such that $\gamma (e)=e.$  Then $m_1=t^k e_{\delta}$ and $m_2=e$ form 
a basis of $\CDcris (D(\delta_1^{-1}))$ and $\Ph (m_2)= m_2-am_1.$ On the other hand, $\Ph$ acts trivially on
 $\CDst (D(\delta_1^{-1}))$ and
we obtain that $a=0$ and $D(\delta_1^{-1})\simeq \CR (\delta)\oplus \CR .$ 

2) Now assume that the proposition holds for $(\Ph,\Gamma)$-modules of rank $d-1.$ Let $D$ be a pure semistable module of rank $d$
with Hodge weights $k_1\leqslant k_2\leqslant \ldots \leqslant k_d.$  Choose a non zero $v\in \F^{k_d}\CDst (D)$ and 
consider the submodule of $D$ which corresponds to $\Bbb Q_p v$ by Proposition 1.2.9:
$$
\CR(\delta_d)=D \cap ( \CR_{\log}[1/t] v).
$$
Then $\delta_d(p)=\lambda p^{-k_d},$ $\delta_d(u)=u^{-k_d}$ ($u\in \Bbb Z_p^*$) and we have an exact sequence
$$
0@>>>\CR(\delta_d)@>>>D@>>>D'@>>>0
$$
where $D'$ is semistable of rank $d-1$ and such that $\CDst (D')^{\Ph=\lambda}=\CDst (D').$ 
 Then $D'\simeq \underset{i=1}\to{\overset{d-1}\to \oplus}
\CR(\delta_i)$ and
$$
\text{Ext}^1(D',\CR (\delta_d))\simeq \underset{i=1}\to{\overset{d-1}\to \oplus}
\text{Ext}^1(\CR (\delta_i)),\CR (\delta_d)).
$$
Let $x$ denote the class of $D$ in $\text{Ext}^1(D',\CR (\delta_d)).$ Since $x$ is semistable,
for any $i$ its image $x_i$ in  $\text{Ext}^1(\CR (\delta_i)),\CR (\delta_d))$ is semistable too.
From  1) it follows that  $x_i=0$ . Thus $x=0$ and $D \simeq D'\oplus \CR (\delta_d).$ The proposition
is proved. 
\enddemo

{\bf 1.5.10.}  Let $D$ be a semistable module.  Assume that  $\CDst (D)=\CDst (D)^{\Ph=1}.$ By Proposition 1.5.9
$D$ is crystalline and
$$
D\simeq \underset{i=1}\to{\overset d \to \oplus }\CR(x^{k_i}),  \qquad k_i\leqslant 0.
$$
In particular $\CDst (D) = D^{\Gamma}$ and the map
$$
i_D\,:\, \CDst (D) \times \CDst (D)@>>> H^1(D) \tag{8}
$$
given by $i_D(\alpha,\beta)=\cl (\alpha,\beta)$ is an isomorphism. 
Now, if  $\CDst (D)^{\Ph=p^{-1}}=\CDst (D)$ then
$\CDst (D^*(\chi))^{\Ph=1}=\CDst (D^*(\chi))$ and we  define $i_D$  by duality.
Let $i_{D,f}$ and $i_{D,c}$ denote the restriction of $i_D$ on the first and the second direct summand 
respectively. Then $\text{Im}(i_{D,f})=H^1_f(D).$ Set $H^1_c(D)=\text{\rm Im} (i_{D,c}).$

\head
{\bf \S2. The $\Cal L$-invariant}
\endhead

{\bf 2.1. Definition of the $\Cal L$-invariant.} 

{\bf 2.1.1.} In this section we generalize Greenberg's definition of the $\Cal L$-invariant. 
Fix a finite set of primes $S$  and denote by ${{\Bbb Q}^{(S)}}/\Bbb Q$  the maximal Galois extension of $\Bbb Q$ 
 unramified outside $S\cup \{\infty\}.$  Set $G_S=\Gal ({\Bbb Q}^{(S)}/\Bbb Q).$
If $M$ is a topological $G_S$-module, we denote by $H^*_S(M)$ the continuous 
cohomology of $G_S$ with coefficients in $M.$
 A $p$-adic representation of $\Gal (\overline{\Bbb Q}/\Bbb Q)$ is said to be 
pseudo-geometric if it satisfies the following conditions:

1) There exists a finite $S$ such that $V$ is unramified outside $S\cup\{\infty\}.$

2) $V$ is potentially semi-stable at $p$.

Let $V$ be a pseudo-geometric representation. Following Bloch and Kato \cite{BK},
for any finite place $v$ we define a subgroup $H^1_f(\Bbb Q_v,V)$ of $H^1(\Bbb Q_v,V)$
by
$$
H^1_f(\Bbb Q_v,V)=\cases 
\ker (H^1(\Bbb Q_v,V)@>>>H^1(\Bbb Q_v^{\text{\rm ur}},V)) & \text{ if $v\ne p$}\\
\ker (H^1(\Bbb Q_v,V)@>>>H^1(\Bbb Q_p,V\otimes_{\Bbb Q_p} \Bc))
& \text{ if $v= p$}.
\endcases
$$
The  Selmer group $H^1_f(V)$ of $V$ is defined as
$$
H^1_f(V)\,=\,\ker \left (
H^1_S(V) @>>>\underset{v\in S}\to \bigoplus 
\frac{H^1(\Bbb Q_v,V)}{H^1_f(\Bbb Q_v,V)}\right ).
$$
Using the inflation-restriction sequence, it is easy to see that this definition does not depend on the choice of $S.$
From the Poitou-Tate exact sequence and the ortogonality of $H^1_f(\Bbb Q_v,V)$ and $H^1_f(\Bbb Q_v,V^*(1))$
one obtains the following exact sequence which relates $H^1_f(V)$ and $H^1_f(V^*(1))$:
$$
\multline
0@>>>H^1_f(V)@>>>H^1_S(V)@>>>\underset{v\in S}\to \bigoplus 
\frac{H^1(\Bbb Q_v,V)}{H^1_f(\Bbb Q_v,V)}@>>>H^1_f(V^*(1))^*@>>>\\
H^2_S(V)@>>>\underset{v\in S}\to \oplus  H^2(\Bbb Q_v,V)
@>>>H^0_S(V^*(1))^*@>>>0 
\endmultline 
\tag{9}
$$
(see \cite{FP}, Proposition 2.2.1). Together with the well known formula
for the Euler characteristic this implies that
$$
\multline
\dim_{\Bbb Q_p} H^1_f(V) - \dim_{\Bbb Q_p} H^1_f(V^*(1)) -\dim_{\Bbb Q_p} H^0_S(V) + \dim_{\Bbb Q_p} H^0_S(V^*(1)) =\\
\dim_{\Bbb Q_p} t_V(\Bbb Q_p)-\dim_{\Bbb Q_p} H^0(\Bbb R,V).
\endmultline
 \tag{10}
$$

\flushpar

{\bf 2.1.2.} Assume that $V$ satisfies the following conditions:

{\bf C1)} $H^1_f(V)=H^1_f(V^*(1))=0$.

{\bf C2)}  $H^0_S(V)=H^0_S(V^*(1))=0$.

{\bf C3)} $V$ is semistable at $p$ and $\Ph \,\,:\,\, \Dst (V)@>>>\Dst (V)$
is semisimple at $1$ and $p^{-1}$.

{\bf C4)} The $(\Ph,\Gamma)$-module $\Ddagrig (V)$ has no crystalline subquotient of the form
$$
0@>>> \Cal R(\vert x\vert x^k)@>>>U @>>>\Cal R@>>>0,\qquad k\geqslant 1.
$$

Set 
$
d_{\pm}(V)=\dim_{\Bbb Q_p}(V^{c=\pm 1}),
$ 
where $c$ denotes the complex conjugation. Comparing {\bf C1-2)} with (10) we obtain that
$
\dim_{\Bbb Q_p} t_V(\Bbb Q_p)=d_{+}(V).
$

{\bf Question 2.1.3.} Let $V$ be an irreducible  pseudo geometric representation which is semistable at $p.$
Does it satisfy {\bf C4)}? This is a straightforward generalization of the hypothesis $\bold U$ 
of \cite{G}. If $V$ is the $p$-adic realization of a pure motive over $\Bbb Q$ of weight $w$ then conjecturally $\Ph \,:\,
\Dc (V)@>>>\Dc (V)$ is semisimple and all of its eigenvalues have the same complex absolute value. In this case
$\bold{C4)}$ holds automatically and $1$ and $1/p$ can not be eigenvalues of $\Ph$ simultaneously.
\newline
\,

 {\bf 2.1.4.} We say that a $(\Ph,N)$-submodule $D$ of $\Dst (V)$ is admissible if the canonical
 projection $D @>>>t_V(\Bbb Q_p)$ is an isomorphism. To any admissible $D$ we associate an increasing filtration
 $(D_i)_{i=-2}^2$ on $\Dst (V)$ by
$$
D_i\,=\,\cases 0 &\text{if $i=-2$,}\\
(1-p^{-1}\Ph^{-1})\,D+N(D^{\Ph=1}) &\text{if $i=-1$,}\\
D &\text{if $i=0$,}\\
D+\Dst (V)^{\Ph=1} \cap N^{-1}(D^{\Ph=p^{-1}})&\text{if $i=1$,}\\
\Dst (V) &\text{if $i=2$}.
\endcases
$$

\proclaim{Lemma 2.1.5} i) $(D_i)_{i=-2}^2$ is the unique filtration on $\Dst (V)$ by $(\Ph,N)$-submodules 
 such that
\newline
\,

{\bf D1)} $D_{-2}=0$, $D_0=D$  and $D_2=\Dst (V);$

{\bf D2)} $(\Dst (V)/D_1)^{\Ph=1,N=0}=0$ and $D_{-1}= (1-p^{-1}\Ph^{-1})D_{-1}\,+N(D_{-1});$

{\bf D3)} $(D_0/D_{-1})^{\Ph=p^{-1}}=D_0/D_{-1}$ and $(D_1/D_0)^{\Ph=1}=D_1/D_0 .$ 
\newline
\,

ii)  Consider the  canonical isomorphism 
$
\Dst (V^*(1))\simeq \text{\rm Hom}_{\Bbb Q_p}(\Dst (V), \Bbb Q_p )
$
where the action of $\Ph$ on the right hand side is given by
$
(\Ph (f))\,(x)\,=\,p^{-1}f(\Ph^{-1}(x)).
$
Set
$$
D^*\,=\, \text{\rm Hom}_{\Bbb Q_p}(\Dst (V)/D, \Bbb Q_p ).
$$
Then $D^*$ is admissible and $D^*_i=\text{\rm Hom}_{\Bbb Q_p}(\Dst (V)/D_{-i}, \Bbb Q_p ).$
\endproclaim

\demo{Proof} i) Since $\Ph$ is semisimple at $1$ and $p^{-1},$ we have a decomposition
of $D$ into a direct sum of $\Ph$-modules
$$
D\simeq X \oplus D^{\Ph=1} \oplus D^{\Ph=p^{-1}}.
$$ 
As $N\Ph=p\,\Ph N,$ one has $N(D^{\Ph=1})\subset D^{\Ph=p^{-1}}$, $N(D^{\Ph=p^{-1}})\subset X$  and $NX\subset X.$
In particular, $D_{-1}=X\oplus D^{\Ph=1}\oplus N(D^{\Ph=1})$ is a $(\Ph,N)$-module 
and $(D/D_{-1})^{\Ph=p^{-1}}=D/D_{-1}.$ Next, as $1-p^{-1}\Ph^{-1}$ is bijective on $X$ one has
$$
(1-p^{-1}\Ph^{-1}) D_{-1}+N(D_{-1})=X+D^{\Ph=1}+ N(D^{\Ph=1})= D_{-1}.
$$
Further, $D_1=X\oplus \left (\Dst (V)^{\Ph=1} \cap N^{-1}(D^{\Ph=p^{-1}}) \right )\oplus  D^{\Ph=p^{-1}}.$ 
From this decomposition it follows immediately that $D_1$ is a $(\Ph,N)$-module such that $(D_1/D_0)^{\Ph=1}= D_1/D_0.$
Now let $N (\bar x)=0$ for some $\bar x=x+D_1\in (\Dst (V)/D_1)^{\Ph=1}.$ As $\Ph$ is semisimple at $1,$
we can assume that $x\in \Dst (V)^{\Ph=1}.$
Then $N x\in D^{\Ph=p^{-1}}$ and we obtain that $x\in \Dst (V)^{\Ph=1} \cap N^{-1}(D^{\Ph=p^{-1}}).$
Thus $\bar x=0$ and $(\Dst (V)/D_1)^{\Ph=1,N=0}=0.$

Conversely, assume that $(D_i)_{i=-2}^2$ is a filtration which satisfies {\bf D1-3)}. From the semisimplicity of $\Ph-1$
it follows that $D_1=D+Y$ where $Y\subset \Dst (V)^{\Ph=1}.$ Since $N(Y)\subset \Dst (V)^{\Ph=p^{-1}} \cap D_1=D^{\Ph=p^{-1}},$
one has $Y\subset \Dst (V)^{\Ph=1} \cap N^{-1}(D^{\Ph=p^{-1}}).$ Let $\bar x=x+D_1$, $x\in \Dst (V)^{\Ph=1}.$
We showed that $N(\bar x)=0$ if and only if $x\in N^{-1}(D^{\Ph=p^{-1}})$ and the condition {\bf D2)} implies that
$Y= \Dst (V)^{\Ph=1} \cap N^{-1}(D^{\Ph=p^{-1}}).$ Thus $D_1=D+ \Dst (V)^{\Ph=1} \cap N^{-1}(D^{\Ph=p^{-1}}).$
A similar argument shows that $D_{-1}=(1-p^{-1}\Ph^{-1})D+N(D^{\Ph=1}).$

ii) The second statement follows from the uniqueness  proved in i) and the fact
that the filtration $\text{\rm Hom}_{\Bbb Q_p}(\Dst (V)/D_{-i}, \Bbb Q_p )$ satisfies {\bf D1-3)}.

\enddemo

{\bf 2.1.6.} Let $D$ be an admissible $(\Ph,N)$-module and $(D_i)_{i=-2}^2$ the associated filtration. By Proposition 1.2.9 it induces 
a filtration of $\Ddagrig (V)$ which we will denote by $(F_i \Ddagrig (V))_{i=2}^2.$ Namely
$$
F_i\Ddagrig (V)\,=\,\Ddagrig (V) \cap (D_i\otimes \Cal R_{\log} \left [1/t \right ]).
$$

Following an idea of Greenberg, define
$$
W=F_1\Ddagrig (V)/F_{-1}\Ddagrig (V).
$$
Then we have an exact sequence
$$
0@>>>\gr_0\Ddagrig (V) @>>>W@>>>\gr_1\Ddagrig (V)@>>>0
\tag{11}
$$
where by Lemma 2.1.5 
$$
\align
&\CDst(\gr_0 \Ddagrig (V))=\CDst(\gr_0 \Ddagrig (V))^{\Ph=p^{-1}},\qquad \F^0\CDst(\gr_0 \Ddagrig (V))=0,\\
&\CDst(\gr_1 \Ddagrig (V))=\CDst(\gr_1 \Ddagrig (V))^{\Ph=1},\qquad \F^0\CDst(\gr_1 \Ddagrig (V))=\CDst(\gr_1 \Ddagrig (V)).
\endalign
$$

\proclaim{Proposition 2.1.7} There exists a unique decomposition
$$
W\simeq W_0\oplus W_1\oplus M
$$
where $W_0$ and  $W_1$ are direct summands of  $\gr_0 \Ddagrig (V)$ and $\gr_1 \Ddagrig (V)$
of ranks $\dim_{\Bbb Q_p} H^0 (W^*(\chi))$ and $\dim_{\Bbb Q_p}H^0(W)$ respectively.
Moreover, $M$ is inserted in an exact sequence
$$
0@>>>M_0@>>>M@>>>M_1@>>>0
$$
where $\gr_0 \Ddagrig (V)\simeq W_0\oplus M_0$,  $\gr_1 \Ddagrig (V)\simeq W_1\oplus M_1$  and
$\text{\rm rg} (M_0)=\text{\rm rg} (M_1).$
\endproclaim
\demo{Proof}  Let $W_1$ denote the saturated $(\Ph,\Gamma)$-submodule of $W$
determined by $\CDst (W_1)=\F^0\CDst (W)^{\Ph=1}$ and let $N_1$ denote the image of
$\CDst (W_1)$ in $\CDst(\gr_1 \Ddagrig (V)).$ From (11) one has  isomorphisms 
\linebreak
$\F^i\CDst (W)\simeq \F^i\CDst (\gr_1 \Ddagrig (V))$ for all $i\geqslant 0.$
This implies that $\CDst (W_1)$ and $N_1$ are isomorphic as filtered Dieudonn\'e modules.
Now remark that in the category of filtered vector spaces every subobject is a direct summand.
As $\Ph$ acts trivially on $\CDst(\gr_1 \Ddagrig (V))$ we obtain that $N_1$ is a direct summand of $\CDst(\gr_1 \Ddagrig (V))$
and there exists a projection of $\CDst (W)$ onto $\CDst (W_1)\simeq N_1.$ This proves that $W_1$ is a direct summand
of $W$. Passing to duals and repeating the same arguments we obtain that $W\simeq W_0\oplus W_1\oplus M.$
It remains to show that  $\text{\rm rg} (M_0)=\text{\rm rg} (M_1).$ As $H^0(M)=0,$ we have an  exact sequence
$$
0@>>>H^0(M_1)@>\delta>>H^1(M_0)@>>>H^1(M).
$$
By  1.5.10 $\dim_{\Bbb Q_p} H^1_f(M_1)=\text{\rm rg} M_1.$
If $\text{\rm rg} (M_0) <\text{\rm rg} (M_1)$ then  there would exist a 
non zero element $\alpha \in H^0 (M_1)$
such that $\delta (\alpha)\in H^1_f(M_0).$ Thus $\delta (\alpha)$ determines a 
crystalline extension of $\CR$ by $M_0$ which is a subquotient of $\Ddagrig (V).$
As $M_0 \simeq \underset{i}\to \oplus \CR (\vert x\vert x^{m_i})$ this violates  {\bf C4)}.
Passing to duals and using the same arguments we obtain the opposite inequality.
The proposition is proved. 
\enddemo 

\proclaim{Lemma 2.1.8} One has 
$$\dim_{\Bbb Q_p}H^1(M)=2r,\quad \dim_{\Bbb Q_p}H^1_f(M)=r,\qquad\text{ where
$r=\text{rg} (M_0)=\text{rg} (M_1).$}
$$
 Moreover
$
H^1_f(M)=\text{\rm Im} (H^1(M_0)@>>>H^1(M)).
$
\endproclaim 
\demo{Proof} By the definition of $W_0$ and $W_1$ we have $H^0(M)=H^0(M^*(\chi))=0$
and the Euler-Poincar\'e characteristic formula gives $\dim_{\Bbb Q_p}H^1(M)=2r.$ 
By Corollary 1.4.5 one has 
$\dim_{\Bbb Q_p}H^1_f(M)=\dim_{\Bbb Q_p}t_M =r.$
Next, as $M$ has no crystalline subquotient of the form {\bf C4)},
 $H^0(M_1)\cap H^1_f(M_0)=\{0\}.$ By 1.5.10 
$$
\dim_{\Bbb Q_p}H^0(M_1)+\dim_{\Bbb Q_p}H^1_f(M_0)=2r =\dim_{\Bbb Q_p}H^1(M_0).
$$
Thus $H^1(M_0)\simeq H^0(M_1)\times H^1_f(M_0)$ and the second statement follows from the fact that
$\dim_{\Bbb Q_p}H^1_f(M_0)=r=\dim_{\Bbb Q_p}H^1_f(M).$ 
\enddemo

{\bf 2.1.9.} Like in \cite{G}, we make the following assumption:
\newline
\,

{\bf C5)} Either $W_0$ or $W_1$ is zero.
\newline
\,

To fix ideas, assume that $W_0=0.$ Otherwise, we  consider $V^*(1)$ instead $V.$ Set $r=\text{rg}(M_0)$,
$s=\text{rg} (W_1)$ and $e=r+s=\text{dim}_{\Bbb Q_p}\left (D_1/D_0\right).$
Consider the short exact sequence
$$
0@>>>F_1\Ddagrig (V)@>>>\Ddagrig (V)@>>>\text{gr}_2\Ddagrig (V)@>>>0
$$

The quotient $\gr_2 \Ddagrig (V)$ is a semistable $(\Ph,\Gamma)$-module
and the admissibility of $D$  implies that all the Hodge weights of  $\CDst (\gr_2 \Ddagrig (V))$ are 
 $\geqslant 0$. Moreover   $\CDst (\gr_2 \Ddagrig (V))^{\Ph=1,N=0}=0$  and by
 Proposition 1.4.4 
$
H^0(\gr_2 \Ddagrig (V))= H^1_f( \gr_2 \Ddagrig (V))=0.
$
By Corollary 1.4.6 one has an exact sequence
$$
0@>>> H^1_f(F_1\Ddagrig (V))@>>>H^1_f(\Ddagrig (V))@>>>H^1_f(\text{gr}_2\Ddagrig (V))=\{0\}
$$
which gives an isomorphism $H^1_f(\Bbb Q_p,V)\simeq H^1_f(F_1\Ddagrig (V)).$
Consider the exact sequence
$$
0@>>>F_{-1}\Ddagrig (V)@>>>F_1\Ddagrig (V)@>>>W@>>>0.
$$
Here $F_{-1}\Ddagrig (V)$ is semistable with Hodge weights $< 0$ and 
$\CDst (F_{-1}\Ddagrig (V))=D_{-1}.$ 
By {\bf D2)}  
\linebreak
$\text{\rm Hom} (D_{-1},\Bbb Q_p)^{\Ph=1,N=0}=0$
and we obtain
$$
\align
&H^0(F_{-1}\Ddagrig (V))=0,\quad H^2(F_{-1}\Ddagrig (V))\simeq H^0 ((F_{-1}\Ddagrig (V))^*(\chi))^*=0,\\
& H^1_f(F_{-1}\Ddagrig (V))=H^1(F_{-1}\Ddagrig (V)).
\endalign
$$
Thus $H^1(W)=\text{coker} (H^1 (F_{-1}\Ddagrig (V))@>>>H^1(F_{1}\Ddagrig (V))).$
On the other hand, by Corollary 1.4.6  $H^1_f(W)=\text{coker} (H^1 (F_{-1}\Ddagrig (V))@>>>H^1_f(F_{1}\Ddagrig (V)))$
and 
$$
\frac{H^1(W)}{H^1_f(W)}\simeq \frac{H^1(F_{1}\Ddagrig (V))}{H^1_f(\Bbb Q_p,V)}.
$$
By Lemma 2.1.8 this is a $\Bbb Q_p$-vector space of dimension $e.$

The conditions {\bf C1-2)} together with the exact sequence (9) give an isomorphism
$$
H^1_S(V)\simeq \underset{v\in S}\to \oplus \, \frac{H^1 (\Bbb Q_v,V)}{H^1_f(\Bbb Q_v,V)}.
$$
Let $H^1(D,V)$ denote the unique subspace of $H^1_S(V)$ whose image under this map is 
\linebreak
$H^1(F_{1}\Ddagrig (V))/H^1_f(\Bbb Q_p,V).$ The localization map 
$H^1_S(V)@>>>H^1 (\Bbb Q_p,V)$ gives rise to a commutative diagram
 $$
 \xymatrix{
H^1(D,V) \ar[d]^{\kappa_D} \ar[dr]^{\bar \kappa_D} & \\
H^1(W) \ar[r] &H^1(\text{gr}_1\Ddagrig (V))}
 $$
 where $\kappa_D$ and $\bar \kappa_D$ are  injective because  by
Lemma 2.1.8  $H^1(\text{gr}_1\Ddagrig (V) ) \simeq {H^1(W)}/{H^1_f(M)}$ and
\linebreak
$H^1_f(M)\subset H^1_f(W).$
Hence, we have a diagram
$$
\xymatrix{
\Cal D_{\text{\rm st}}({\text{\rm gr}}_1 \Ddagrig (V)) \ar[r]^{\overset{i_{D,f}}\to \sim} & H^1_f({\text{\rm gr}}_1 \Ddagrig (V))\\
H^1(D,V) \ar[u]^{\rho_{D,f}} \ar[r]^{\bar \kappa_D} \ar[d]_{\rho_{D,c}} &
H^1({\text{\rm gr}}_1 \Ddagrig (V)) \ar[u]_{p_{D,f}}
\ar[d]^{p_{D,c}}
\\
\Cal D_{\text{\rm st}}({\text{\rm gr}}_1 \Ddagrig (V)) \ar[r]^{\overset{i_{D,c}}\to\sim} &H^1_c({\text{\rm gr}}_1 \Ddagrig (V)),}
$$
where $\rho_{D,f}$ and $\rho_{D,c}$ are defined as the unique maps making
this diagram commute.
From the definition of $H^1(D,V)$ it follows that  $\rho_{D,c}$ is an isomorphism.

\proclaim{Definition} The determinant
$$
\Cal L (V,D)= \det \left ( \rho_{D,f} \circ \rho^{-1}_{D,c}\,\mid \,\Cal
D_{\text{\rm st}}(\text{\rm gr}_1 \Ddagrig (V)) \right )
$$
will be called  the $\Cal L$-invariant associated to $V$ and $D$.
\endproclaim

The isomorphism $\gr_1\Ddagrig (V)\simeq M_1\oplus W_1$ induces a decomposition
$$
H^1 (D,V)\simeq H^1_{\text{loc}} (D,V) \oplus H^1_{\text{gl}} (D,V)
$$
such that $\rho_{D,c}(H^1_{\text{loc}} (D,V))= \CDst (M_1)$ and 
$\rho_{D,c}(H^1_{\text{gl}} (D,V))= \CDst (W_1)$.

\proclaim{Proposition 2.1.10} One has
$$
\bar \kappa_D (H^1_{\text{loc}}(D,V)) =\text{\rm Im} (H^1(M)@>>>H^1(M_1)).
$$

\endproclaim
\demo{Proof} Consider the diagram
$$
\xymatrix{
0 \ar[r] & H^0(M_1)\ar[r] & H^1(M_0) \ar[r] & H^1(M) \ar[r] & H^1(M_1) \ar[r] & H^2(M_0) \ar[r] & 0\\
& & & H^1_{\text{loc}}(D,V)\ar[u]^{\kappa_D} \ar[ur]^{\bar \kappa_D} &  & &
}
$$ 
As $H^0(M)=H^2(M)=0,$ the upper row is exact. This implies that the image of $H^1(M)$ in $H^1(M_1)$
is a $\Bbb Q_p$-vector space of dimension $r.$ On the other hand, $\dim  H^1_{\text{loc}}(D,V)=r$
and the proposition follows from the commutativity of the diagram.
\enddemo

\proclaim{Corollary 2.1.11} One has
$$
\Cal L(V,D)=\Cal L_{\text{\rm loc}}(V,D)\, \Cal L_{\text{\rm gl}}(V,D)
$$
where $\Cal L_{\text{\rm loc}}(V,D)=\det \left ( \rho_{D,f} \circ \rho^{-1}_{D,c}\,\mid \,\Cal
D_{\text{\rm st}}(M_1) \right )$ and  $\Cal L_{\text{\rm gl}}(V,D)=\det \left ( \rho_{D,f} \circ \rho^{-1}_{D,c}\,\mid \,\Cal
D_{\text{\rm st}}(W_1) \right ).$ Here the factor  $\Cal L_{\text{\rm loc}}(V,D)$ is local, i.e. depends only 
on the restriction of the representation $V$ on a decomposition group at $p$.
\endproclaim

{\bf 2.2. $p$-adic $L$-functions.}

{\bf 2.2.1.} Let $V$ be a pseudo-geometric representation which satisfies {\bf C1-5)}. It seems reasonable
to expect that to any admissible subspace $D$ of $\Dst (V)$ 
one can associate a meromorphic $p$-adic $L$-function $L^{\text{an}}_{p}(V,D,s)$ interpolating special values
of the complex $L$-function $L(V,s).$ 
A natural generalization of Greenberg's conjecture
about the behavior of $L_p^{\text{an}}(V,D,s)$  states as follows. If $f$ is a semisimple endomorphism of a vector space $W$
then $W\simeq \ker (f)\oplus \text{Im} (f)$ and we define ${\det}^* (f)\,=\,\det (f\,\vert \, \text{Im} (f)).$
\newline
\,

{\bf Greenberg's conjecture.} $L_{p}(V,D,s)$ has a zero of order $e$ at $s=0$.
 Moreover, $\Cal L(V,D)\ne 0$ and
$$
\lim_{s\to 0} \frac{L_{p}^{\text{an}}(V,D,s)}{s^{e}}\,= \,\Cal L (V,D)\,\Cal E(V,D)  \,\frac{L(V,0)}{\Omega_{\infty}(V)}
$$
where 
$$\Cal E(V,D)= 
{\det}^* \left (1-p^{-1}\Ph^{-1} \,\vert D^{N=0} \right )\,  {\det}^* \left ( 1-p^{-1}\Ph^{-1} \,\vert \,(D^*)^{N=0} \right )$$
and  $\Omega_{\infty}(V)$ denotes the  Deligne's period.
\newline
\,

{\bf Remarks 2.2.2.} 1) The interpolation factor $\Cal E(V,D)$ can be written in the form
$$
\Cal E(V,D)={\det}^* \left (1-p^{-1}\Ph^{-1} \,\vert D^{N=0} \right )\,\det 
\left (1-\Ph\,\vert \, \Dst (V)/(N\Dst (V)+D)\right ).
$$
Remark that $\Cal E(V,D)=\Cal E(V^*(1),D^*).$

2) Assume that $V$ is ordinary at $p$. Then $V$ is equipped with an increasing filtration
$\Cal F_iV$ such that $\text{gr}_i (V)\,(-i)$ are unramified. Set $D= \Dst (\Cal F_{-1}V).$ Then  $\Cal L(V,D)$
coincides with  Greenberg's $\Cal L$-invariant and the above conjecture coincides 
with  the conjecture formulated  in \cite{G}, p. 166. 

3) Our definition of the $\Cal L$-invariant generalizes without modifications to the case of an arbitrary 
coefficient field $L/\Bbb Q_p.$
\newline
\,

{\bf 2.2.3.}  Let $V$ be the $p$-adic representation associated to a normalized newform $f$ of weight $2k$ 
on $\Gamma_0(Np)$ ($(N,p)=1$) with Fourier coefficients in $L/\Bbb Q_p$. Assume that $U_p(f)=p^{k-1}f$ where $U_p$ is the Atkin-Lehner operator.
Then $V$ is semistable of Hodge-Tate weights $(0,2k-1)$ and the $(\Ph,N)$-module $\Dst (V)$ is generated
by two elements $d_1$ and $d_2$ such that $\Ph (d_2)=p^k d_2,$ $\Ph (d_1)=p^{k-1}d_1$, $N d_2=d_1$,
$N d_1=0.$ The $\Cal L$-invariant of Fontaine-Mazur is defined as the unique element $\Cal L_{\text{FM}} (f)\in  L$
such that
$
d_2-\Cal L_{\text{FM}} (f)\, d_1\in \F^0 \Dst (V).
$

Then $\Dst (V(k))$ is generated by $d^*_1$ and $d^*_2$ such that $\Ph (d^*_2)= d_2^*,$
$\Ph (d^*_1)=p^{-1}d^*_1$, $N(d^*_2)=d^*_1$ and $D=L d^*_1$ is the unique admissible  $(\Ph,N)$-submodule
of $\Dst (V(k)).$ 

\proclaim{Proposition 2.2.4} One has  $\Cal L_{\text{FM}} (f) =\Cal L(V(k),D).$
\endproclaim
\demo{Proof} 
By Proposition 1.3.2 the refinement $\Cal F_0 \Dst (V) =L d_1$, $\Cal F^1 \Dst (V)=\Dst (V)$
gives rise to an exact sequence
$$
0@>>>\CR_{L}(\vert x\vert^{-(k-1)})@>>>\Ddagrig (V)@>>>\CR_{ L}(x^{-(2k-1)}\vert x\vert ^{-k})@>>>0.
$$
We have canonical isomorphisms
$$
\Ext^1_{\CR_{L}}(\CR_{L}(x^{-(2k-1)}\vert x\vert ^{-k}),\CR_{L}(\vert x\vert^{-(k-1)}))\simeq 
H^1(\CR_{L} (x^{2k-1}\vert x\vert )) \simeq
H^1(C^\bullet_{\text{st}}(\CR_{L} (x^{2k-1}\vert x\vert ))).
$$
By the proof of Proposition 1.4.4 the  class of $\Dst (V)$ in 
$H^1(C^\bullet_{\text{st}}(\CR_{L} (x^{2k-1}\vert x\vert )))$ is
\linebreak
 $\cl (-\Cal L_{\text{FM}} (f),0,1).$
On the other hand, we have an exact sequence
$$
0@>>>\CR_{L}(x^k\vert x\vert)@>>>\Ddagrig (V(k))@>>>\CR_{ L}(x^{1-k})@>>>0.
$$
which coincides with the sequence (11) for $V(k).$ Denote by $X$ the image of
$H^1(\Ddagrig (V(k)))$ in 
\linebreak
$H^1(\CR_{L}(x^{-(k-1)})).$ Then
$\Cal L(D,V(k))$ is the unique element  $\lambda \in L$ such that $y^*_{k-1}-\lambda x^*_{k-1} \in X.$
Let $\delta^*_{k-1}\,:\,H^0(\CR_{ L}(x^{-(k-1)}))@>>>H^1(\CR_{L}(x^k\vert x\vert))$ 
denote the connecting map associated to the dual exact sequence
$$
0@>>>\CR_{L}(x^k\vert x\vert)@>>>\Ddagrig (V^*(1-k))@>>>\CR_{ L}(x^{1-k})@>>>0.
$$
By local duality and Proposition 1.5.7 i) $\text{Im} (\delta^*_{k-1})$ is generated by
$\beta^*_k-\lambda \alpha^*_k.$ Consider the commutative diagram
$$
\xymatrix{
0 \ar[r] &\CR_L(x^{2k-1}\vert x \vert) \ar[d]^{t^{k-1}} \ar[r] &\Ddagrig (V^*(1-k)) (x^{k-1}) \ar[r] \ar[d]^{t^{k-1}}
& \CR_L \ar[r] \ar[d]^{t^{k-1}} &0\\
0 \ar[r] &\CR_L(x^{k}\vert x \vert)  \ar[r] &\Ddagrig (V^*(1-k))  \ar[r] 
& \CR_L (x^{1-k})\ar[r]  &0}
$$
where the vertical maps are the multiplication by $t^{k-1}.$
This induces a commutative diagram
$$
\xymatrix{
H^0(\CR_L) \ar[r]^{\delta_0^*} \ar[d]^{t^{k-1}} &H^1(\CR_L(x^{2k-1}\vert x\vert))\ar[d]^{t^{k-1}}\\
H^0(\CR_L(x^{1-k})) \ar[r]^{\delta^*_{k-1}} &H^1(\CR_L(x^k\vert x\vert)).}
$$
Together with Corollary 1.5.6 this shows that $\text{Im} (\delta_0^*)$ is generated by $\beta^*_{2k-1}-\lambda \alpha^*_{2k-1}.$ 
On the other hand  $\text{Im} (\delta_0^*)$ is generated by the class of $\Ddagrig (V^*(1-k))$ in 
$$
\Ext^1_{\CR_{L}}(\CR_{L}(x^{-(k-1)}),\CR_{L}(x^k\vert x\vert))\simeq 
H^1(\CR_{L} (x^{2k-1}\vert x\vert )).
$$
A  short direct computation gives $\cl (\Ddagrig (V^*(1-k)))=
-\cl (\Ddagrig (V(k))).$ But by Proposition 1.5.8 ii) we have 
$$
\cl (\Ddagrig (V(k)))=-\beta_{2k-1}^*+\Cal L_{\text{FM}}(f) \,\alpha^*_{2k-1}.
$$
This proves that $\lambda =\Cal L_{\text{FM}}(f).$

\enddemo

{\bf 2.2.5.} Assume that $V$ is crystalline at $p$.
In  this  case Greenberg's conjecture is compatible with Perrin-Riou's theory of 
$p$-adic $L$-functions. Namely,  set $\Gamma_1=\text{Gal}(\Bbb Q_p(\zeta_{p^{\infty}})/\Bbb Q_p(\zeta_p))$ and fix
a generator $\gamma_1\in \Gamma_1$. 
Let $\Lambda =\Bbb Z_p[[\Gamma_1]]$ denote the Iwasawa algebra of $\Gamma_1$ and 
$\Cal H(\Gamma_1)=\{f(\gamma_1-1)\,\vert \, f\in \Cal H\}$ where
$\Cal H$ is  the algebra of power series which converges on the open unit disc.
Denote by $\Cal K (\Gamma_1)$ the field of fractions of $\Cal H .$
Let $D$ be a $\Ph$-stable  subspace of 
 $\Dc (V)$ such that $D\oplus \F^0 \Dc (V)=\Dc (V).$ In \cite{PR}, Perrin-Riou constructed a free $\Lambda$-submodule
$\bold L_{\Iw}(V,D)$ of $\Cal K(\Gamma_1)$ in terms of Galois cohomology of $V$.
This construction depends on the choice of a $G_S$-stable lattice  $T$ of $V$ and a $\Bbb Z_p$-lattice $N$ of $D$
but for simplicity we do not take it into account here and assume that the $p$-adic period associated to this choice is
equal to $1$.
 The main conjecture
of the Iwasawa theory of $V$ states as follow. 
\newline
\newline
{\bf Main Conjecture} (\cite{PR}, \cite{Cz1}). One has $L_p^{\text{an}}(V,D,s)= f(\chi (\gamma_1)^s-1)$ where $f$ is an 
appropriate generator of $\bold L_{\Iw}(V,D)$. 
\newline
\,

To fix ideas assume that $D^{\Ph=1}=0.$  Then $e= \dim_{\Bbb Q_p}(D^{\Ph=p^{-1}}).$ In \cite{Ben2} we prove the following result.

\proclaim{Theorem 2.2.6} Assume that $\Cal L(V,D)\ne 0.$ Let $f$ be a generator of $\bold L_{\Iw}(V,D).$ Then
  $f(\chi (\gamma_1)^s-1)$ is a meromorphic $p$-adic function which has a zero 
at $s=0$ of order $e$ and
$$
\lim_{s\to 0} \frac{f(\chi (\gamma_1)^s-1) }{s^{e}}    \overset{p}\to\sim 
 \,\Cal L (V,D)\, \Cal E(V,D)\,  \frac{\#\sha (T)}
{\# H^0_S(V/T) \,\# H^0_S(V^*(1)/T^*(1))}\,\Tam^0 {(T)}\,
$$
where $\sha (T)$ is the Tate-Shafarevich group of Bloch and Kato (see \cite{FP}, section 5.3.4)  
and $\Tam^0(T)$ is the product of local Tamagawa numbers.
\endproclaim

Remark that the Bloch-Kato conjecture predicts that 
$$
\frac{L(V,0)}{\Omega_{\infty}(V)}\,=\,\frac{\#\sha (T)}
{\# H^0_S(V/T) \,\# H^0_S(V^*(1)/T^*(1))}\,\Tam^0 {(T)}
$$
and  Theorem 2.2.6 implies the compatibility of the Greenberg conjecture with
Perrin-Riou's theory.

\enddocument
\bye

This is a remark about maximal crystalline/semistable submodules.

\proclaim{Lemma 1.3.8} Let $D$ be a $(\Ph,\Gamma)$-module over $\CR$. Then there exists  unique
saturated submodules $D^{\text{cr}}$ and $D^{\text{st}}$ of $D$ such that

1) $D^{\text{cr}}$ is crystalline and $D^{\text{st}}$ is semistable;

2) If $D'$ is a crystalline (resp. semistable) submodule of $D$, then
$D'\subset D^{\text{cr}}$ (resp. $D'\subset D^{\text{st}}$).

Moreover $\CDcris (D^{\text{cr}})=\CDcris (D)$ and $\CDst (D^{\text{st}})=\CDst (D).$
\endproclaim 
\demo{Proof} Let $D_1$ and $D_2$ be two semistable submodules of $D$. Then $\widetilde D=D_1+D_2$ is a
$(\Ph,\Gamma)$-module which seats in the exact sequence
$$
0@>>>D_1 \cap D_2@>>> D_1\oplus D_2 @>>>\widetilde D @>>>0.
$$
By Lemma 1.3.5, $\widetilde D$ is semistable. Let $\widetilde D^{\text{sat}}$ denote the saturation
of $\widetilde D$ in $D$. Then $\text{rg} ( \widetilde D^{\text{sat}})=\text{rg} (\widetilde D).$
On the other hand, by \cite{BC}, Lemma ??? $  \widetilde D^{\text{sat}}=\widetilde D[1/t] \cap D.$
Thus $\CDst (\widetilde D^{\text{sat}})=\CDst (\widetilde D)$ and we obtain that
$\widetilde D^{\text{sat}}$ is a saturated semistable submodule of $D$ containing $D_1$ and $D_2$.
This proves the existence of $D^{\text{st}}\subset D$ satisfying 1) and 2). 
Next, the inclusion $\CDst (D)\otimes_K \CR_{\log}[1/t] \hookrightarrow D \otimes_{\CR}\CR_{\log} [1/t]$
gives
$$
\bold D[1/t] \,=\,(\CDst (D)\otimes_K \CR_{\log}[1/t])^{N=0}  \hookrightarrow D[1/t]
$$ 
and from the construction of $\Cal M$ (see the proof of Proposition 1.3.7) it follows that $\Cal M(\CDst (D))\subset D.$
Thus $\Cal M (\CDst (D))= D^{\text{st}}$ and $\CDst (D^{\text{st}})=\CDst (D).$ 

The existence and properties of $D^{\text{cr}}$ can be proved by the same arguments. 
\enddemo

\proclaim{Corollary 1.3.9} The functor $D\mapsto D^{\text{cr}}$ (resp. $D\mapsto D^{\text{st}}$)
is the right-adjoint of the embedding $j_{\text{\rm cris}}\,:\,{\bold M}^{\,\Ph,\Gamma}_{\text{\rm cris},K}@>>>
{\bold M}^{\,\Ph,\Gamma}_{K}$ (resp. 
$j_{\text{\rm st}}\,:\,{\bold M}^{\,\Ph,\Gamma}_{\text{\rm st},K}@>>>
{\bold M}^{\,\Ph,\Gamma}_{K}$). 
\endproclaim
\demo{Proof} Let $D_1$ and $D_2$ be two $(\Ph,\Gamma)$-modules. Assume that $D_1$ is crystalline.
If $f\in \Hom (j_{\text{cris}}(D_1), D_2)$, then $f(D_1)\subset D_2$ is crystalline.
Thus, $f(D_1)\subset D_2^{\text{cr}}$ and we obtain that
$$
 \Hom (j_{\text{cris}}(D_1), D_2)=
 \Hom(D_1, D_2^{\text{st}})
$$
 In the semistable case the proof is analogous.
\enddemo

\proclaim{Lemma 1.4.2} Let $D$ be a $(\Ph,\Gamma)$-module. Then 

i) $\cl (a,b)\in H^1_f(D)$ (resp. $\cl (a,b)\in H^1_{st}(D)$) if and only if
the equation $(\gamma-1)\,x=b$ has a solution in $D[1/t]$ (resp. in $D\otimes_{\CR}\CR_{\log} [1/t]$).

ii) The natural maps $H^1_{f}(D^{\text{cr}})@>>>H^1_f(D)$ and $H^1_{st}(D^{\text{st}})@>>>H^1_{st}(D)$ are
injective.
\endproclaim
\demo{Proof} i) An extension $D_\alpha$ is crystalline (resp. semistable) if and only if
there exists $x\in D[1/t]$ (resp. $x\in D\otimes_{\CR}\CR_{\log} [1/t]$) such that 
$x+e \in \CDcris (D_\alpha)$ (resp. $x+e\in \CDst (D_\alpha)$). As $(\gamma -1)\,e=b,$ this proves i).
The second statement can be easily deduced from i) but we will give another proof
which can be applied to a more general situation. 
Let $\cl (\alpha)\in \text{\rm Ext}^1 (\CR,D^{\text{cr}})$ and let $\cl (\beta)$ denote
the image of $\cl(\alpha)$ in $\text{\rm Ext}^1 (\CR,D).$ Then we have a commutative
diagram
$$
\xymatrix{
0 \ar[r] &D^{\text{cr}} \ar[r] \ar[d]^{i} &D_{\alpha}  \ar[r] \ar[d] \ar @{.>}[dl]^{s}&\CR \ar[r] \ar[d]^{=}&0\\
0\ar[r]&D \ar[r] &D_{\beta}\ar[r]  &\CR \ar[r] & 0.
}
$$
Assume that $\cl (\alpha)$ is crystalline and  $\cl (\beta)=0.$
Then the bottom exact sequence splits and there exists a map $s\,:\,D_\alpha @>>> D$
such that $s \mid_{D^{\text{cr}}} =i.$ Since the $(\Ph,\Gamma)$-module $D_\alpha$ is a crystalline,  by Corollary 1.3.9
$s(D_\alpha)\subset D^{\text{cr}}$ 
and $s$ defines a splitting of the top row. This proves that the map
$H^1_f(D^{\text{cr}})@>>> H^1_f(D)$ is injective. 
The semistable version can be proved by the same arguments.

\enddemo

\enddocument
\bye